\documentclass[11pt,reqno]{amsart}
\usepackage[letterpaper,margin=1.05in]{geometry}
\usepackage[T1]{fontenc}
\usepackage{lmodern,microtype,mathtools,amssymb,amsthm,mathrsfs}
\usepackage{enumitem,xcolor,hyperref}
\definecolor{linkblue}{RGB}{24,55,90}
\hypersetup{colorlinks=true,linkcolor=linkblue,citecolor=linkblue,
urlcolor=linkblue,pdftitle={Holder estimates for complex Hessian equations on Kahler manifolds}}
\setlist[enumerate]{leftmargin=*,itemsep=4pt}
\numberwithin{equation}{section}
\newtheorem{theorem}{Theorem}[section]
\newtheorem{proposition}[theorem]{Proposition}
\newtheorem{lemma}[theorem]{Lemma}
\newtheorem{corollary}[theorem]{Corollary}
\theoremstyle{definition}
\newtheorem{definition}[theorem]{Definition}
\theoremstyle{remark}
\newtheorem{remark}[theorem]{Remark}
\newcommand{\ddc}{i\partial\bar\partial}
\newcommand{\R}{\mathbb R}

\newcommand{\Sk}{S_k}
\DeclareMathOperator{\tr}{tr}
\DeclareMathOperator{\osc}{osc}
\allowdisplaybreaks[2]
\title[H\"older estimates on K\"ahler manifolds]
{Global $C^{2,\alpha}$ Estimates for Complex \(k\)-Hessian\\
Equations on K\"ahler Manifolds}
\author{Yulun Xu}
\date{}
\begin{document}
\begin{abstract}
We prove an a priori \(C^{2,\alpha}\) estimate for classical admissible
solutions of the complex \(k\)-Hessian equation on a compact K\"ahler
manifold with a positive 
\(C^\alpha\) right-hand side.
\end{abstract}
\maketitle

\setcounter{tocdepth}{1}
\tableofcontents
\clearpage

\section{Introduction}\label{sec:statement}
Throughout this paper, $(M,\chi)$ denotes a connected compact K\"ahler manifold of complex dimension $n \ge 2$.  Let \(\omega\) be a fixed smooth
real \emph{closed} \((1,1)\)-form with
\(\lambda_\chi(\omega)\in\Gamma_k\), where
\[
 \Gamma_k=\{\lambda\in\R^n:\sigma_j(\lambda)>0,\ 1\le j\le k\},
 \qquad 1\le k\le n.
\]
For a real function \(u\), put
$
 \omega_u=\omega+\ddc u.
$
In this paper, we study the complex k-Hessian equation:
\begin{equation}\label{eq:main}
 \omega_u^k\wedge\chi^{n-k}=e^G\chi^n,\qquad
 \lambda_\chi(\omega_u)\in\Gamma_k,\qquad
 G\in C^\alpha(M),\qquad \sup_Mu=0.
\end{equation}

The complex \(k\)-Hessian equations interpolate between the Poisson
equation and the complex Monge--Amp\`ere equation, which correspond to the cases $k=1$ and $k=n$ respectively. 

For the complex Monge--Amp\`ere equation, Yau's solution of the Calabi
conjecture established smooth solvability for smooth positive
densities satisfying the volume constraint \cite{Yau1978}. Calabi \cite{Ca} developed the third derivative estimate for real Monge-Amp\'ere equations, and Yau \cite{Yau1978} and Riebesehl-Schulz \cite{RS} established its complex counterparts. Siu \cite{siu} developed the complex-geometric form of the Evans-Krylov argument. For less regular densities, Kolodziej \cite{ko} proved a uniform estimate for the potential when $e^G \in L^p$ for $p>1$. Chen-He \cite{CH} proved gradient and Laplacian estimates for $G \in W^{1,p}, p >2n$. Chen-Cheng 
\cite{CC} proved a gradient estimate under a square-Dini condition on the modulus of continuity of $F$. Wu-Zheng \cite{WuZheng} proved a gradient estimate under the assumption that $e^{\frac{G}{\beta}} \in W^{1,\beta}$ for some $\beta>n$. 

For complex k-Hessian equations with general $k$, on compact K\"ahler manifolds, Hou--Ma--Wu obtained
a second-order estimate controlled by the gradient
\cite{HMW}. Dinew--Ko{\l}odziej supplied the gradient estimate
through a Liouville theorem and a blow-up argument, completing
smooth solvability in the standard K\"ahler setting
\cite{DK2017}. Recently, Chu-Liu \cite{ChuLiu} obtained another proof of gradient estimate without using blow-up argument. Sz\'ekelyhidi  developed a general
theory of fully nonlinear elliptic equations including complex $k$-Hessian equations on compact Hermitian
manifolds based on \(\mathcal C\)-subsolutions \cite{Szekelyhidi2018}. Zhang \cite{Zhang} independently solved complex $k$-Hessian equations on compact Hermitian manifolds in the case $\chi=\omega$. 

In the case that the density function has lower regularity, Blocki \cite{Blocki} introduced weak solutions. Uniform estimate, stability, and weak existence were proved by Dinew-Kolodziej \cite{DK2014}.  The author developed quantitative
regularization of admissible functions and proved H\"older
continuity of pluripotential solutions of complex \(k\)-Hessian
equations on compact Hermitian manifolds for densities in
\(L^p\), \(p>n/k\) \cite{Xu2026}.  More recently, Guo-Song \cite{GuoSong} proved gradient estimate depending on $||G||_{C^1}$ without using Hou-Ma-Wu estimate.

The recent theorem of Li--Wang--Zheng gives a global
\(C^{2,\alpha}\) estimate for the complex Monge--Amp\`ere
equation with positive \(C^\alpha\) density on compact K\"ahler
manifolds \cite{LWZ2026}. 

Parabolic methods also form part of this background.
Sheng--Wang studied a complex Hessian flow on compact K\"ahler
manifolds and proved long-time existence and convergence
\cite{ShengWang2019}. Phong--T\^o developed a general parabolic
\(\mathcal C\)-subsolution framework on compact Hermitian
manifolds \cite{PhongTo2021}. 

In this paper, we generalize Li-Wang-Zheng's global $C^{2,\alpha}$ estimate for Monge-Amp\'ere equation \cite{LWZ2026} to complex $k$-Hessian equations:

\begin{theorem}\label{thm:main}
Let $(M,\chi)$ be a K\"ahler manifold and let \(0<\alpha<1\). For any $G \in C^{\alpha}(M)$ with $\int_M e^G \chi^n = \int_M \omega^k \wedge \chi^{n-k}$,  there exists a unique solution 
\(u\in C^{2,\alpha}(M,\R)\) satisfying (\ref{eq:main}), and
\begin{equation}\label{eq:target}
 \|u\|_{C^{2,\alpha}(M)}
 \le C(M,\chi,\omega,k,\alpha,\|G\|_{C^\alpha(M)}).
\end{equation}
\end{theorem}

The proof proceeds in four steps:

\begin{enumerate}
\item Using the linearized operator of $\log \sigma_k$, we can define a Hermitian metric $a$. Similar to Guo-Phong-Song-Sturm's Sobolev inequality for Monge-Amp\'ere equations \cite{GuoPhongSongSturm}, we use auxiliary function to prove a Sobolev estimate for $a$.
\item We use Guo-Song's idea of using a stability estimate to prove a $C^1$ estimate for complex $k$-Hessian equations \cite{GuoSong}. Using the Sobolev inequality that we establish, we can apply the integral method, and this allows us to get $C^1$ estimate using $||G||_{C^{\alpha}}$, while the $C^1$ estimate by Guo-Song \cite{GuoSong} depends on $||e^G||_{C^1}$, where $e^G$ may degenerate. 
\item We use the idea of unshifted restart in Li-Wang-Zheng \cite{LWZ2026} which proves $C^{1,\beta}$ estimate for $\beta$ close to $1$. Roughly speaking, we define the original flow:
\begin{equation}
 \psi_t= \log \frac{\omega_{\psi}^k \wedge \chi^{n-k}}{\chi^n}-P_tG,\qquad \psi(0)=u,
\end{equation}
and its unshifted restart at an original time \(\tau>0\):
\begin{equation}
 \partial_s\varphi^\tau=\log \frac{\omega_{\varphi^\tau}^k \wedge \chi^{n-k}}{\chi^n} -P_sG,
 \qquad \varphi^\tau(0)=\psi(\tau).
\end{equation}
Using the $C^1$ estimate on $u$, we can derive $C^1$ estimate on $\psi$ and $\varphi^\tau$. We can also get 
\begin{equation}\label{e initial}
||\varphi^{\tau}-\psi||_{\infty}\le C \tau^{1+\alpha/2}.
\end{equation}

Using spatial interpolation with parabolic Moser iteration, we can obtain:
\begin{equation}\label{e tr est}
||tr_{\chi}\omega_{\psi}||_{\infty}(t)\le C_{1,\epsilon} t^{-d_{1}-\epsilon}, \,\,\, d_{1}= 1,
\end{equation}
for any $\epsilon>0$.

In order to have a better estimate on $tr_{\chi}\omega+ \Delta_\chi \psi$, we can use the good estimate (\ref{e initial}) to improve the estimate between two unshifted restart. Then, the exponent in (\ref{e tr est}) can be improved to be $- d_{j+1} -\epsilon$, where $d_{j+1}= d_{j} / (1+ \alpha/2)$. Doing induction, we can get that:
\begin{equation}
||tr_{\chi}\omega_{\psi}||_{\infty}(t) \le C_{\epsilon} t^{-\epsilon},\end{equation}
for any $\epsilon>0$. This implies that $u \in C^{1,\beta}$ for any $\beta<1$.

\item Li-Wang-Zheng \cite{LWZ2026} used Li-Li-Zhang's results \cite{LLZ2018} to improve the $C^{1,\beta}$ estimate to $C^{2,\alpha}$ estimate. We use again the flows above to generalize Li-Li-Zhang's result from the complex Monge-Amp\'ere equation to the complex $k$-Hessian equation. 

We obtain the \(W^{2,p}\) estimate for \(u\) by combining the trace bound in
Proposition~\ref{lem:gradient-bootstrap} with the volume estimate for
\(\{\tr A_u>2\tr A_\psi(t)\}\) in Lemma~\ref{lem:initial-trace-superlevel}.
For constant scalar curvature K\"ahler metrics, Chen--Cheng obtain a pointwise
trace bound from a finite trace moment using the scalar curvature equation
\cite[Proposition~4.2]{ChenCheng2021}. For the complex \(k\)-Hessian
equation, we compare \(\psi(s)\) with the fixed potential \(\psi(2t)\) on
\([t/2,t]\). The parabolic maximum principle bounds the maximum trace at
\(t\) by its maximum at \(2t\) plus \(Ct^{\alpha/12}\).
Summing \eqref{eq:trace-two-times} over dyadic times gives the uniform
bound for \(\tr_\chi\omega_u\) in Proposition~\ref{prop:trace-cap}.

\end{enumerate}

In Section 2, we introduce some basic tools and definitions. In Section 3, we prove Sobolev inequality for a metric defined using complex $k$-Hessian equation. In section 4, we prove $C^1$ estimate for complex $k$
-Hessian equation. In Section 5, we prove a spatial $C^{1,\beta}$ estimate of flow for any $\beta<1$. In Section 6, we prove $W^{2,p}$ estimate for complex $k$-Hessian equation. In Section 7, we prove $W^{2,p}$ estimate and trace propagation along flow. In Section 8, we prove a uniform estimate on $tr_{\chi} \omega_u$, where $u$ is the solution to the complex $k$-Hessian equation. In Section 9, we prove Theorem \ref{thm:main}.

\paragraph{\textbf{Declaration on the use of AI.}}
ChatGPT 6.1 Sol and ChatGPT 6.2 Astra were used as an adversarial reader and to improve the writing of the arguments. In particular, ChatGPT provided a simple way to prove $W^{2,p}$ estimate and the uniform complex Hessian estimate.\section{Preliminaries}

\subsection{Notations}
For any function $v\in C^2(M)$, we denote $A_v = \chi^{-1}\omega_v$.  For any Hermitian matrix A, we denote $\Sk(A)=\binom nk^{-1}\sigma_k(A).$

We use \(\Delta=\Delta_\chi=\chi^{i\bar j}\partial_i\partial_{\bar j}\)
and \(\partial_i=(\partial_{x_i}-i\partial_{y_i})/2\).

\subsection{Properties of $\sigma_k$}
We recall the generalized Newton--Maclaurin inequality; see
\cite{Spruck2005}. We use the convention \(\sigma_0=1\).
\begin{proposition}\label{prop:newton-maclaurin}
For $\lambda\in\Gamma_k$ and $k>l\geq0$, $r>s\geq0$,
$k\geq r$, $l\geq s$, we have
\begin{equation}\label{eq:newton-maclaurin}
\left[
\frac{\sigma_k(\lambda)/\binom nk}
     {\sigma_l(\lambda)/\binom nl}
\right]^{\frac{1}{k-l}}
\leq
\left[
\frac{\sigma_r(\lambda)/\binom nr}
     {\sigma_s(\lambda)/\binom ns}
\right]^{\frac{1}{r-s}}.
\end{equation}
If \((k,l)\ne(r,s)\), equality holds if and only if
\(\lambda_1=\lambda_2=\cdots=\lambda_n>0\).
If \((k,l)=(r,s)\), equality is automatic.
\end{proposition}

\begin{definition}[Uniform comparability]\label{def:uniform-comparability}
For nonnegative quantities \(Q_1,Q_2\), we write \(Q_1\asymp Q_2\)
if there are constants \(0<c\le C<\infty\) such that
\[
 cQ_2\le Q_1\le CQ_2
\]
throughout the domain under consideration. Unless otherwise stated,
these comparison constants depend only on
\(M,\chi,\omega,n,k,\alpha,\|G\|_{C^\alpha(M)}\).
Any additional dependence is specified where the comparison is used.
The constants are uniform in the spatial and time variables and in
the choice of smooth approximation.
\end{definition}

\subsection{Heat flow}
Let \(P_t=e^{t\Delta_\chi}\) be the heat semigroup of the fixed
K\"ahler metric \(\chi\). We first record its regularization estimates.

\begin{lemma}\label{lem:heat-regularization}
Let \(0<\alpha<1\) and \(G\in C^\alpha(M)\).
The semigroup \(P_t\) is positive, preserves constants and
\(\chi^n\), contracts the supremum norm, and is self-adjoint
in \(L^2(\chi^n)\). For \(0<t\le1\),
\begin{equation}\label{eq:heat-envelope}
 \|P_tG-G\|_\infty\le Ct^{\alpha/2},\quad
 \|\partial P_tG\|_\infty\le Ct^{(\alpha-1)/2},\quad
 \|\nabla\bar\nabla P_tG\|_\infty+\|\partial_tP_tG\|_\infty
 \le Ct^{\alpha/2-1}.
\end{equation}
Here \(C=C(M,\chi,n,\alpha)\|G\|_{C^\alpha(M)}\);
no higher derivative of \(G\) is required.
\end{lemma}
\begin{proof}
Let \(H_t(x,y)\) be the heat kernel relative to \(dV=\chi^n\).
The maximum principle gives positivity and supremum contraction.
The symmetry of \(\Delta_\chi\), together with \(\Delta_\chi1=0\),
gives self-adjointness, constant preservation, and volume preservation.
In particular,
\[
 P_tG(x)-G(x)=\int_M H_t(x,y)\bigl(G(y)-G(x)\bigr)\,dV(y).
\]
The Gaussian heat-kernel bounds on the fixed compact manifold imply
\[
\begin{aligned}
 \int_M H_t(x,y)d_\chi(x,y)^\alpha\,dV(y)&\le Ct^{\alpha/2},\\
 \int_M\bigl(\sqrt t\,|\partial_xH_t|
       +t|\nabla_x\bar\nabla_xH_t|+t|\partial_tH_t|\bigr)
       d_\chi(x,y)^\alpha\,dV(y)&\le Ct^{\alpha/2}.
\end{aligned}
\]
For the differentiated kernels, the integral against \(1\) is zero,
so \(G(y)\) can again be replaced by \(G(y)-G(x)\).
Applying \(|G(y)-G(x)|\le C\|G\|_{C^\alpha}d_\chi(x,y)^\alpha\)
to these formulas proves \eqref{eq:heat-envelope}.
\end{proof}

\begin{lemma}\label{lem:heat-c1beta-laplacian}
For \(0<\beta<1\), \(h\in C^{1,\beta}(M)\), and \(0<\eta\le1\),
\begin{equation}\label{eq:heat-c1beta-laplacian}
 \|\Delta_\chi P_\eta h\|_\infty
 \le C_\beta\eta^{-(1-\beta)/2}\|h\|_{C^{1,\beta}(M)},
\end{equation}
where \(C_\beta=C(M,\chi,n,\beta)\).
\end{lemma}
\begin{proof}
Fix \(x_0\in M\) and real normal coordinates
\(z=(z^1,\ldots,z^{2n})\) centered at \(x_0\) on
\(B_{2r}(x_0)\), where \(r>0\) is fixed below the injectivity radius.
Choose a smooth cutoff \(\zeta_{x_0}\), equal to one on
\(B_r(x_0)\) and supported in \(B_{2r}(x_0)\), and define
\[
 \ell_{x_0}(y)=h(x_0)+\zeta_{x_0}(y)
 \sum_{a=1}^{2n}\frac{\partial h}{\partial z^a}(x_0)
       \bigl(z^a(y)-z^a(x_0)\bigr).
\]
The cutoff term is extended by zero outside the coordinate ball.
Throughout the proof \(x_0\) is fixed and all differentiations act
on \(y\), so the coefficients
\(\partial h/\partial z^a(x_0)\) are constants. In particular,
\[
 \Delta_\chi\ell_{x_0}
 =\sum_{a=1}^{2n}\frac{\partial h}{\partial z^a}(x_0)
   \Delta_\chi\!\left[
   \zeta_{x_0}\bigl(z^a-z^a(x_0)\bigr)\right],
 \qquad
 \|\Delta_\chi\ell_{x_0}\|_\infty\le C\|h\|_{C^1(M)}.
\]
The constants are uniform in \(x_0\), since the metric and radius
are fixed and \(M\) is compact; no second derivative of \(h\) occurs.
Taylor's formula in these coordinates gives
\[
 |h(y)-\ell_{x_0}(y)|
 \le C_\beta\|h\|_{C^{1,\beta}(M)}
          d_\chi(x_0,y)^{1+\beta}.
\]
This holds first on \(B_r(x_0)\); outside that ball it follows
from the uniform bound on \(h-\ell_{x_0}\) and
\(d_\chi(x_0,y)\ge r\).

The Gaussian heat-kernel derivative bounds on the fixed compact
manifold, using \cite[Corollary~1.2]{Hsu1999}, give
\[
 |\Delta_{\chi,x}H_\eta(x,y)|
 \le C\eta^{-n-1}
      \exp\!\left(-\frac{d_\chi(x,y)^2}{C\eta}\right).
\]
Since \(P_\eta\) commutes with \(\Delta_\chi\) on smooth functions,
\[
\begin{aligned}
 \Delta_\chi P_\eta h(x_0)
 &=\int_M\Delta_{\chi,x}H_\eta(x_0,y)
       \bigl(h(y)-\ell_{x_0}(y)\bigr)\,dV(y)\\
 &\quad+P_\eta(\Delta_\chi\ell_{x_0})(x_0).
\end{aligned}
\]
Here the derivative acts only on the kernel, not on the base-point
dependence of \(\ell_{x_0}\). The Gaussian moment estimate
\[
 \int_M e^{-d_\chi(x_0,y)^2/(C\eta)}
            d_\chi(x_0,y)^{1+\beta}\,dV(y)
 \le C_\beta\eta^{n+(1+\beta)/2}
\]
follows by the change of variables \(z=\sqrt\eta\,w\) in
\(B_r(x_0)\); the complement contributes an exponentially small
term. Supremum contraction of \(P_\eta\) therefore yields
\[
\begin{aligned}
 |\Delta_\chi P_\eta h(x_0)|
 &\le C_\beta\|h\|_{C^{1,\beta}(M)}
          \eta^{-n-1}\eta^{n+(1+\beta)/2}+C\|h\|_{C^1(M)}\\
 &\le C_\beta\eta^{-(1-\beta)/2}\|h\|_{C^{1,\beta}(M)}.
\end{aligned}
\]
The last inequality uses \(0<\eta\le1\). Taking the supremum over
\(x_0\) proves \eqref{eq:heat-c1beta-laplacian}.
\end{proof}

\subsection{Estimates related to linearized operator}\label{sec:algebra}
For a \(\chi\)-Hermitian endomorphism \(\mathsf H\), we denote by
\(|\mathsf H|_{\mathrm F}\) its Frobenius norm induced by \(\chi\).
In a \(\chi\)-unitary frame, \(\mathsf H\) is represented by a Hermitian
matrix and
\[
 |\mathsf H|_{\mathrm F}
 =\left(\sum_{i,j}|\mathsf H_{i\bar j}|^2\right)^{1/2}
 =\bigl(\tr(\mathsf H^*\mathsf H)\bigr)^{1/2}
 =\left(\sum_{i=1}^n\lambda_i(\mathsf H)^2\right)^{1/2},
\]
where \(\mathsf H^*\) is the adjoint with respect to \(\chi\).

For an admissible real \((1,1)\)-form \(\Omega\), write
\(\mathsf H=\chi^{-1}\Omega\), set
\(F(\mathsf H)=\log \sigma_k(\mathsf H)\). The first variation
\[
 \left.\frac{d}{dt}F\bigl(\chi^{-1}(\Omega+tH)\bigr)\right|_{t=0}
=\sum_{i,j}F^{i\bar j}H_{i\bar j}
\]
for every real \((1,1)\)-form \(H\), with \(\chi\) fixed. 

Define
\[
Lh=F^{i\bar j}h_{i\bar j},\quad
 \Gamma(h,\ell)=\Re(F^{i\bar j}h_i\ell_{\bar j}),\quad
 \Gamma(h)=\Gamma(h,h),\quad \Gamma_{\chi}(h,l)= \Re(\chi^{i\bar j}h_i\ell_{\bar j}).
\]

In a \(\chi\)-unitary frame diagonalizing \(\mathsf H\), $(F^{i\bar{j}})$ is diagonal and
\[
 F^{i \bar i}=\frac{\partial F}{\partial\lambda_i}
     =\frac{\sigma_{k-1}(\lambda(\mathsf H)|i)}{\sigma_k(\mathsf H)}.
\]
Here \(\lambda(\mathsf H)|i\) denotes the eigenvalue vector with its
\(i\)-th entry omitted.  Define
\begin{align*}
 \mathcal F&:=\tr((F^{i\bar{j}}))=\sum_i F^{i \bar i}.
\end{align*}

\begin{lemma}\label{lem:algebra}
Let \(2\le k\le n\), \(\lambda(\mathsf H)\in\Gamma_k\). Then
\begin{enumerate}[label=\textup{(\arabic*)},ref=\arabic*]
\item\label{item:algebra-norm} \(|\mathsf H|_{\mathrm F}<\tr \mathsf H\).
\item\label{item:algebra-ellipticity} \(F^{i \bar i}\ge(\tr \mathsf H)^{-1}\) for every \(i\).
\item\label{item:algebra-quadratic}
\(\displaystyle\sum_i F^{i \bar i}\lambda_i^2\ge(\tr \mathsf H)/n\).
\end{enumerate}
\end{lemma}
\begin{proof}
Since \(\Gamma_k\subset\Gamma_2\),
\((\tr \mathsf H)^2-|\mathsf H|_{\mathrm F}^2=2\sigma_2(\mathsf H)>0\) and \(\tr \mathsf H>0\).
This proves the first assertion.

For the second, the cone deletion property gives
\(\lambda|i\in\Gamma_{k-1}\), and
\[
 F^{i \bar i}\tr \mathsf H-1
 =\frac{\sigma_1(\lambda|i)\sigma_{k-1}(\lambda|i)
                  -\sigma_k(\lambda|i)}{\sigma_k(\lambda)}.
\]
The numerator is positive if \(\sigma_k(\lambda|i)\le0\).
Otherwise \(k<n\) and \(\lambda|i\in\Gamma_k\). Applying
Proposition~\ref{prop:newton-maclaurin} in dimension \(n-1\) gives
\[
 \frac{\sigma_k(\lambda|i)}{\sigma_{k-1}(\lambda|i)}
 \le\frac{n-k}{k(n-1)}\sigma_1(\lambda|i)
 \le\sigma_1(\lambda|i).
\]
Thus \(F^{i \bar i}\tr \mathsf H\ge1\).

Finally, expanding \(\sigma_1\sigma_k\) and dividing by \(\sigma_k\)
yields
\[
 \sum_i F^{i \bar i}\lambda_i^2
 =\tr \mathsf H-(k+1)\frac{\sigma_{k+1}(\lambda)}{\sigma_k(\lambda)}.
\]
If \(\sigma_{k+1}\le0\), the right-hand side is at least \(\tr \mathsf H\).
If \(\sigma_{k+1}>0\), then \(k<n\) and \(\lambda\in\Gamma_{k+1}\);
Proposition~\ref{prop:newton-maclaurin} gives
\[
 (k+1)\frac{\sigma_{k+1}(\lambda)}{\sigma_k(\lambda)}
 \le\frac{n-k}{n}\tr \mathsf H.
\]
For \(k=n\), use \(\sigma_{n+1}=0\). These cases prove the third
assertion, and in particular \(\sum_i F^{i \bar i}\lambda_i^2\ge(\tr \mathsf H)/n\).
\end{proof}

\begin{lemma}\label{lem:trace-growth}
Let \(2\le k\le n\) and \(\lambda(\mathsf H)\in\Gamma_k\).
On a density window \(0<f_-\le\Sk(\mathsf H)\le f_+\),
\begin{equation}\label{eq:trace-growth}
 \mathcal F\le\frac{k}{f_-}(\tr \mathsf H)^{k-1},
\end{equation}
and
\begin{equation}\label{eq:dual-growth}
 \tr \mathsf H\le\frac{n}{k^{k-1}}f_+(\mathcal F)^{k-1}.
\end{equation}
\end{lemma}
\begin{proof}
Each monomial of degree \(k-1\) occurs in \(n-k+1\) deleted
symmetric functions. Hence
\[
 \mathcal F=\frac{(n-k+1)\sigma_{k-1}(\mathsf H)}{\sigma_k(\mathsf H)}.
\]
Lemma~\ref{lem:algebra} (1) gives \(|\lambda_i|<\tr \mathsf H\), which proves
\eqref{eq:trace-growth}. For the reverse estimate, put
\(b_j=\sigma_j(\mathsf H)/\binom nj\). For \(k>2\),
Proposition~\ref{prop:newton-maclaurin} gives
\(b_{k-1}^{k-1}\ge b_1b_k^{k-2}\); for \(k=2\), this is an identity. Since
\(\mathcal F=kb_{k-1}/b_k\),
\[
 (\mathcal F)^{k-1}\ge\frac{k^{k-1}b_1}{b_k}
 =\frac{k^{k-1}\tr \mathsf H}{n\Sk(\mathsf H)},
\]
which proves \eqref{eq:dual-growth}.
\end{proof}

\begin{remark}
Although Lemma \ref{lem:algebra} and Lemma \ref{lem:trace-growth} assume that $k \ge 2$, it doesn't cause any loss of generality. In fact, Theorem \ref{thm:main} holds trivially in the case $k=1$ which is the Poisson equation. As a result, we can assume that $k \ge 2$ throughout this paper.
\end{remark}

\begin{lemma}\label{lem:determinant}
For \(\lambda(\mathsf H)\in\Gamma_k\),
\begin{equation}\label{eq:determinant}
 det(F^{i\bar{j}}) \ge(k/n)^n\Sk(\mathsf H)^{-n/k}.
\end{equation}
\end{lemma}
\begin{proof}
For every positive definite Hermitian \(P\),
concavity, homogeneity, and Proposition~\ref{prop:newton-maclaurin}
give
\[
 \tr(\nabla_{\mathsf H}(\Sk^{1/k})(\mathsf H)P) =  \Sk(\mathsf H)^{1/k} + \tr(\nabla_{\mathsf H}(\Sk^{1/k})(\mathsf H)(P- \mathsf H)\ge \Sk(P)^{1/k}\ge(\det P)^{1/n}.
\]
Since \((F^{i \bar j})=k\nabla_{\mathsf H}(\Sk^{1/k})(\mathsf H)/\Sk(\mathsf H)^{1/k}\), the matrix arithmetic--geometric
mean inequality (Proposition~\ref{prop:newton-maclaurin}) gives
\[
 n(\det (F^{i \bar j}))^{1/n}=\min_{P>0,\,\det P=1}\tr((\nabla F) P)\ge k/\Sk(\mathsf H)^{1/k}.
\]
Raising the last inequality to the \(n\)-th power proves
\eqref{eq:determinant}.
\end{proof}

\begin{lemma}\label{lem:reverse-background}
There is \(\delta_\omega>0\), depending only on \(M,\chi,\omega,n,k\),
such that every admissible \(u\) satisfies
\begin{equation}\label{eq:potential}
 F^{i \bar j} \omega_{i \bar j} \ge\delta_\omega\mathcal F,
 \qquad Lu=k- F^{i \bar j} \omega_{i \bar j}
             \le k-\delta_\omega\mathcal F.
\end{equation}
Here \(F^{i \bar j}\) and \(L\) are evaluated at \(A_u\).
\end{lemma}
\begin{proof}
Strict admissibility and compactness give \(\delta_\omega>0\)
such that \(\chi^{-1}(\omega-\delta_\omega\chi)\in\Gamma_k\).
Concavity and homogeneity of the root operator imply
\[
 \tr\bigl(\nabla_{\mathsf H}(\Sk^{1/k})(A_u)\chi^{-1}(\omega-\delta_\omega\chi)\bigr)
 \ge \Sk\bigl(\chi^{-1}(\omega-\delta_\omega\chi)\bigr)^{1/k}>0.
\]
Multiply by \(k/\Sk(A_u)^{1/k}\) and use Euler's identity
\(\tr(aA_u)=k\) to obtain \eqref{eq:potential}.
\end{proof}

\begin{lemma}\label{lem:newton-symmetry}
Let \(\Omega\) be a smooth closed strictly admissible form, and put
\(\mathsf H=\chi^{-1}\Omega\), \(f=\Sk(\mathsf H)\), \(d\nu=f\chi^n\).
Here \(a,L,\Gamma\) are evaluated at this \(\mathsf H\). Then
\begin{equation}\label{eq:newton-forms}
 (Lh)\,d\nu=k\,\ddc h\wedge\Omega^{k-1}\wedge\chi^{n-k},
\end{equation}
and, for smooth real \(h,\ell\),
\begin{equation}\label{eq:ibp}
 \int_MhL\ell\,d\nu=-\int_M\Gamma(h,\ell)\,d\nu,
 \qquad\int_ML\ell\,d\nu=0.
\end{equation}
The fixed Laplacian satisfies
\begin{equation}\label{eq:fixed-ibp}
 \int_Mh\Delta\ell\,dV=-\int_M\Gamma_\chi(h,\ell)\,dV,
 \qquad dV=\chi^n.
\end{equation}
\end{lemma}
\begin{proof}
Differentiate \(f\chi^n=\Omega^k\wedge\chi^{n-k}\) in the
direction \(\ddc h\) to obtain \eqref{eq:newton-forms}.
The form \(\Omega^{k-1}\wedge\chi^{n-k}\) is closed, so integration
by parts gives the first identity in \eqref{eq:ibp}. Take \(h=1\)
for the second. The same argument with \(\Omega=\chi\), \(k=1\),
gives \eqref{eq:fixed-ibp}.
\end{proof}

\section{Sobolev inequality for the complex $k$-Hessian equation}\label{sec:sobolev}
We first control the trace introduced by a partition of unity.

\begin{lemma}\label{lem:trace-form}
Let \(a\) be smooth and positive definite, let \(d\nu\) be a smooth
positive measure, and suppose \(L\) satisfies \eqref{eq:ibp}.
If a smooth real \(u\) satisfies
\[
 \osc_Mu\le K,\qquad Lu\le B-\delta\mathcal F,\qquad\delta>0,
\]
then every smooth real \(z\) satisfies
\begin{equation}\label{eq:trace-form}
 \int(\mathcal F)z^2\,d\nu
 \le C\left(\int\Gamma(z)\,d\nu+\|z\|_{L^2(\nu)}^2\right),
\end{equation}
where \(C=C(K,\max(B,0),\delta)\).
\end{lemma}
\begin{proof}
Set \(h(u)=2-e^{-(u-\inf u)}\), so \(1\le h<2\),
\(h'\ge e^{-K}\), and \(h''=-h'\). Therefore
\[
 \frac{-Lh}{h}=\frac{h'}h[-Lu+\Gamma(u)]
 \ge\frac{\delta e^{-K}}2\mathcal F-\max(B,0).
\]
Using \eqref{eq:ibp} and completing the square gives
\[
 \int z^2\frac{-Lh}{h}\,d\nu
 =\int\Gamma(z)\,d\nu-\int|dz-z\,d\log h|_a^2\,d\nu
 \le\int\Gamma(z)\,d\nu.
\]
Combining the two inequalities proves \eqref{eq:trace-form}.
\end{proof}

\begin{lemma}\label{lem:local-sobolev}
Let \(a\) be smooth and positive definite. Suppose \(L\) is symmetric as in \eqref{eq:ibp},
\(d\nu=f\chi^n\), \(0<f_-\le f\le f_+\), and
\(\det_\chi a\ge d>0\). On a fixed coordinate ball \(D\)
compactly contained in a holomorphic chart,
\begin{equation}\label{eq:local-sobolev}
 \|z\|_{L^p(\nu)}^2\le C\int_D\Gamma(z)\,d\nu,
 \qquad z\in C_c^\infty(D,\R),\quad2<p<\frac{2n}{n-1}.
\end{equation}
For \(n=1\), any finite \(p>2\) is allowed. The constant depends
only on \(p,f_-,f_+,d\) and the fixed chart and metric.
\end{lemma}
\begin{proof}
Choose \(r>n\), smooth \(g\ge0\) with \(\|g\|_{L^r(\nu)}=1\),
and \(0<\epsilon\le1\). The determinant assumption gives
\(\det(a^{i\bar j})\ge d/\sup_{\overline D}\det(\chi_{i\bar j})\).
Solve
\[
 \det(\psi_{i\bar j})=\frac{(g+\epsilon)^n\sup_{\overline D}\det(\chi_{i\bar j})}{n^nd},
 \qquad\psi|_{\partial D}=0.
\]
The Dirichlet estimate at exponent \(r/n>1\), with homogeneous
\(1/n\) dependence \cite[Theorem 1.2]{WangWangZhou2021}, gives \(-C\le\psi\le0\).
The matrix arithmetic--geometric mean inequality, a special case of
Proposition~\ref{prop:newton-maclaurin}, gives \(L\psi\ge g+\epsilon\).
Thus \(H=1-\psi\) satisfies \(1\le H\le H_0\) and \(-LH\ge g\),
with uniform \(H_0\). By \eqref{eq:ibp},
\[
 H_0^{-1}\int gz^2\,d\nu
 \le\int z^2\frac{-LH}{H}\,d\nu
 =\int\Gamma(z)\,d\nu-\int|dz-z\,d\log H|_a^2\,d\nu
 \le\int\Gamma(z)\,d\nu.
\]
Taking the supremum over \(g\), using density in \(L^r\), proves
\eqref{eq:local-sobolev} for \(p=2r/(r-1)\). Varying \(r>n\)
gives the stated range.
\end{proof}

\begin{proposition}\label{prop:sobolev}
Assume the hypotheses of Lemma~\ref{lem:trace-form},
\(d\nu=f\chi^n\) with \(0<f_-\le f\le f_+\), and
\(\det_\chi a\ge d>0\). Then, for \(2<p<2n/(n-1)\),
\begin{equation}\label{eq:sobolev}
 \|z\|_{L^p(\nu)}^2\le C\left(\int\Gamma(z)\,d\nu
                                  +\|z\|_{L^2(\nu)}^2\right).
\end{equation}
For \(n=1\), every finite \(p>2\) is allowed. The constant depends only on
\[
 M,\chi,p,f_-,f_+,d,K,\max(B,0),\delta.
\]
It uses no upper or
derivative bound for \(a\), and no derivative bound for \(f\).
\end{proposition}
\begin{proof}
Choose smooth real cutoffs \(\eta_\ell\), subordinate to a finite
coordinate cover, with \(\sum_\ell\eta_\ell^2=1\).
The mixed terms in \(\sum_\ell\Gamma(\eta_\ell z)\) cancel because
\(\sum_\ell\eta_\ell d\eta_\ell=0\). The cutoff gradients are fixed,
so
\[
 \sum_\ell\int\Gamma(\eta_\ell z)\,d\nu
 \le\int\Gamma(z)\,d\nu+C\int(\mathcal F)z^2\,d\nu.
\]
The triangle inequality in \(L^{p/2}\),
Lemma~\ref{lem:local-sobolev}, and Lemma~\ref{lem:trace-form} give
\[
 \|z\|_{L^p(\nu)}^2\le\sum_\ell\|\eta_\ell z\|_{L^p(\nu)}^2
 \le C\left(\int\Gamma(z)\,d\nu+\int z^2\,d\nu\right).
\]
This is \eqref{eq:sobolev}.
\end{proof}

\begin{remark}
For complex Monge-Amp\'ere equation, the Sobolev inequality was proved by Guo-Phong-Song-Sturm \cite{GuoPhongSongSturm}, whose Sobolev constant depends only on the $L^p$ norm of the right hand side of the equation.
\end{remark}

\begin{proposition}\label{prop:hessian-sobolev}
Let \((M,\chi)\) be a compact K\"ahler manifold without boundary,
fix \(1\le k\le n\), and let \(\omega\) be a smooth closed real \((1,1)\)-form with
\(\lambda_\chi(\omega)\in\Gamma_k\). Suppose \(u\in C^\infty(M,\R)\)
satisfies
\[
 \lambda_\chi(\omega_u)\in\Gamma_k,\qquad
 \osc_Mu\le K,\qquad
 0<f_-\le f:=\Sk(A_u)\le f_+.
\]
Set \(a=\nabla_{\mathsf H}\log\Sk(A_u)\) and \(d\nu=f\chi^n\).
Then \eqref{eq:sobolev} holds for every smooth real \(z\) and
\(2<p<2n/(n-1)\), with
\(C=C(M,\chi,\omega,k,p,K,f_-,f_+)\).
For \(n=1\), every finite \(p>2\) is allowed.
The constant requires no derivative bound for \(f\) and no
upper or derivative bound for \(a\).
\end{proposition}
\begin{proof}
Admissibility makes \(a\) smooth and positive definite.
Since \(\omega_u\) and \(\chi\) are closed,
Lemma~\ref{lem:newton-symmetry} gives the symmetry identity
\eqref{eq:ibp} with respect to \(d\nu=f\chi^n\).
Lemma~\ref{lem:determinant} and \eqref{eq:potential} give
\[
 \det_\chi a\ge(k/n)^nf_+^{-n/k}>0,\qquad
 Lu\le k-\delta_\omega\mathcal F,
\]
where \(\delta_\omega>0\) depends only on \(M,\chi,\omega,n,k\).
Thus \(a\), \(d\nu\), and \(u\) satisfy every hypothesis of
Proposition~\ref{prop:sobolev}: symmetry, the density bounds,
the determinant lower bound, and the oscillation and potential
bounds, with \(B=k\) and \(\delta=\delta_\omega\).
That proposition proves the claim with the stated dependencies.
\end{proof}

\section{$C^1$ estimate for complex $k$-Hessian equation}\label{sec:gradient}
In the elliptic estimates of this section, use
\(A_u=\chi^{-1}(\omega+\ddc u)\), and evaluate \(a,L,\Gamma\)
at \(A_u\). The parabolic case is specified separately in
Corollary~\ref{cor:gradient-flow}.
We first obtain potential compactness from Xu's H\"older estimate.

\begin{lemma}\label{lem:potential-holder}
Every smooth normalized solution of \eqref{eq:main} satisfies
\[
 \|u\|_{C^{1/(n+1)}(M)}
 \le C(M,\chi,\omega,n,k,\|G\|_{C^0(M)}).
\]
Consequently, a family of these solutions with a common
\(C^0\) bound for \(G\) is precompact in \(C^0(M)\).
In particular, \(-C_0\le u\le0\) uniformly.
\end{lemma}
\begin{proof}
Apply the estimate in \cite[Theorem 1.1]{Xu2026} to
\(\Sk^{1/k}\) on \(\Gamma_k\). This operator is positive,
concave, homogeneous of degree one, elliptic on \(\Gamma_k\),
and zero on its boundary. Since \(\nabla_{\mathsf H}(\Sk^{1/k})=(\Sk^{1/k}/k)a\),
Lemma~\ref{lem:determinant} gives the required structural bound
\[
 \det_\chi(\nabla_{\mathsf H}(\Sk^{1/k}))
 =(\Sk^{1/k}/k)^n\det_\chi a\ge n^{-n}.
\]
The fixed background is strictly admissible. The smooth function
\(u-1\), normalized by \(\sup_M(u-1)=-1\), is a viscosity
solution of the root equation \(\Sk(A_{u-1})^{1/k}=e^{G/k}\).
In Xu's notation choose \(p_0=2\), so \(q_0=2\).
The norm \(\|e^{nG/k}\|_{L^2(\chi^n)}\) is controlled by
\(M,\chi,n,k\) and \(\|G\|_{C^0}\), and
\(1/(n+1)<2/(1+2n)\). Xu's theorem therefore gives the
stated H\"older bound. Arzel\`a--Ascoli proves the compactness
assertion, and \(\sup_Mu=0\) gives the potential bounds.
\end{proof}

For the comparison argument, choose \(a_0,\beta>0\) in that order.
Given a smooth admissible \(w\) with \(\|u-w\|_\infty\le\beta^{-1}\),
write, within this section,
\begin{equation}\label{eq:gradient-weight}
\begin{split}
 q&=|\partial u|_\chi^2+K_0,\qquad
 \phi=-a_0u+\beta(u-w-\beta^{-1})^2,\qquad W=e^\phi q \\
 K_0&=K_{a_0}(1+\|\partial w\|_\infty)^4.
 \end{split}
\end{equation}
Here \(K_{a_0}\) is a scalar constant chosen in
Lemma~\ref{lem:gradient-bochner}.
These quantities are used only in the gradient estimates below.

\begin{lemma}\label{lem:kahler-bochner}
Let \(u\) be smooth and real, and let
\(Lh=a^{i\bar j}h_{i\bar j}\), with \(a\) positive definite.
At any point choose K\"ahler normal coordinates of $\chi$ in which
\(a^{i\bar j}=F^{i \bar i}\delta_{ij}\), and write
\(u_{pi}=\nabla_i u_p\), \(u_{p\bar i}=\nabla_{\bar i}u_p\).
Use the curvature convention
\(R_{i\bar j p\bar q}=-\partial_i\partial_{\bar j}\chi_{p\bar q}\)
at the center of these coordinates. Then
\begin{equation}\label{eq:kahler-bochner}
\begin{aligned}
 L|\partial u|_\chi^2
 ={}&\sum_{i,p}F^{i \bar i}\bigl(|u_{pi}|^2+|u_{p\bar i}|^2\bigr)
       +2\Re\sum_{i,p}F^{i \bar i} u_{\bar p}\nabla_pu_{i\bar i}\\
 &+\sum_{i,p,q}F^{i \bar i}R_{i\bar i p\bar q}u_q u_{\bar p}.
\end{aligned}
\end{equation}
\end{lemma}
\begin{proof}
Differentiate
\(|\partial u|_\chi^2=\chi^{p\bar q}u_pu_{\bar q}\)
twice in the chosen coordinates. The first derivatives of \(\chi\)
vanish there. The products of second derivatives of \(u\) give
the two squared-Hessian terms. The third ordinary derivatives give
the real part in \eqref{eq:kahler-bochner}, since
\(\nabla_pu_{i\bar i}=\partial_p\partial_i\partial_{\bar i}u\)
at the chosen point. The second derivatives of
\(\chi^{-1}\) give the displayed curvature contraction.
\end{proof}

\begin{lemma}\label{lem:gradient-bochner}
Let \(u\) be smooth and strictly admissible, with \(-C_0\le u\le0\) and
\(0<f_-\le\Sk(A_u)\le f_+\).
There are choices of \(a_0,\beta,K_{a_0},c,C>0\), depending only on
\(M,\chi,\omega,n,k,C_0,f_-,f_+\), such that every smooth admissible
\(w\) with \(\|u-w\|_\infty\le\beta^{-1}\), and the associated
weights \eqref{eq:gradient-weight}, satisfy
\begin{equation}\label{eq:gradient-comparison}
 LW\ge ce^\phi\tr A_u+cW\{\mathcal F+\Gamma(u-w)\}-CW
            +2e^\phi\Re\langle\partial u,\partial\log f\rangle_\chi,
 \qquad f=\sigma_k(A_u).
\end{equation}
\end{lemma}
\begin{proof}
Work in a K\"ahler normal frame diagonalizing \(A_u\). Write
\[
 E_{A_u}=\sum_i F^{i \bar i}\lambda_i^2,\qquad
 E_u=\sum_{i,p}F^{i \bar i}|u_{p\bar i}|^2.
\]
We will choose \(K_{a_0}\ge1\), so \(q\ge1+|\partial u|_\chi^2\).
The comparison radius gives
\(-2\le\beta(u-w-\beta^{-1})\le0\). Differentiating \(\phi\) yields
\[
 \phi_i=-a_0w_i+
 \bigl[-a_0+2\beta(u-w-\beta^{-1})\bigr](u-w)_i.
\]
By \eqref{eq:potential} and
\(L(u-w)=k-\tr(aA_w)\le k\), its second derivative satisfies
\begin{equation}\label{eq:gradient-weight-laplacian}
 L\phi\ge(a_0\delta_\omega)\mathcal F-k(a_0+4)
                 +2\beta\Gamma(u-w).
\end{equation}

Differentiating the equation gives
\[
 \sum_i F^{i \bar i}\nabla_pu_{i\bar i}
 =(\log f)_p-\sum_i F^{i \bar i}\nabla_p\omega_{i\bar i}.
\]
Substitute this identity into \eqref{eq:kahler-bochner}.
The curvature and background-form terms have absolute value at most
\[
 C_n\bigl(\|\mathrm{Rm}(\chi)\|_\infty|\partial u|_\chi^2
       +\|\nabla\omega\|_\infty|\partial u|_\chi\bigr)\mathcal F
 \le Cq\mathcal F.
\]
Here \(C\) depends only on
\(n,\|\mathrm{Rm}(\chi)\|_\infty,\|\nabla\omega\|_\infty\).
Thus
\begin{equation}\label{eq:bochner}
 L|\partial u|_\chi^2\ge\sum_{i,p}F^{i \bar i}|u_{pi}|^2+E_u
 +2\Re\langle\partial u,\partial\log f\rangle_\chi-Cq\mathcal F.
\end{equation}
Also \(u_{p\bar i}=\lambda_i\delta_{pi}-\omega_{p\bar i}\)
implies \(E_u\ge E_{A_u}/2-C_\omega\mathcal F\), where \(C_\omega\)
depends only on \(n,\|\omega\|_{L^\infty(\chi)}\).
The product rule for \(W=e^\phi q\), followed by completing the
square in the pure Hessian, now gives
\begin{equation}\label{eq:gradient-completed-square}
\begin{aligned}
 e^{-\phi}LW\ge{}&
 \sum_{i,p}F^{i \bar i}|u_{pi}+u_p\phi_i|^2
 +E_u+qL\phi+K_0\Gamma(\phi)\\
 &+2\Re\sum_i F^{i \bar i}\lambda_i u_i\phi_{\bar i}
 -2\Re\sum_{i,p}F^{i \bar i}\omega_{i\bar p}u_p\phi_{\bar i}\\
 &-Cq\mathcal F
 +2\Re\langle\partial u,\partial\log f\rangle_\chi.
\end{aligned}
\end{equation}
Indeed the completed square uses
\(\sum_p|u_p|^2=q-K_0\), leaving \(K_0\Gamma(\phi)\).

Apply Lemma~\ref{lem:gradient-mixed} below to the two mixed terms.
The first uses \(E_{A_u}/8\); the second uses \(q\Gamma(u-w)/4\).
Since \(|\partial u|_\chi^2+(1+\|\partial w\|_\infty)^2\le q\), their remaining
\(\Gamma(u-w)\) contribution is at most \(C_{a_0}q\Gamma(u-w)\).
Choose \(a_0\) so that \(a_0\delta_\omega\) exceeds the fixed
Bochner coefficient of \(q\mathcal F\) by \(4\). Next choose
\(\beta\) so that \(2\beta-C_{a_0}-1/4\ge1\).
These choices are independent of derivatives of \(w\).

Finally choose \(K_{a_0}\) large. The inequalities
\[
\begin{gathered}
 \frac{(1+\|\partial w\|_\infty)^4}{q}\le K_{a_0}^{-1},\qquad
 \frac1q\le K_{a_0}^{-1},\\
 \frac{(1+\|\partial w\|_\infty)|\partial u|_\chi}{q}
 \le\frac1{2\sqrt{K_{a_0}}(1+\|\partial w\|_\infty)}
 \le\frac1{2\sqrt{K_{a_0}}}
\end{gathered}
\]
make the remaining mixed-term errors and
\(C_\omega\mathcal F\) smaller than the reserved multiple of
\(q\mathcal F\). Discarding the nonnegative squares and
\(K_0\Gamma(\phi)\), we obtain
\[
 e^{-\phi}LW\ge
 \tfrac14E_{A_u}+q\{\mathcal F+\Gamma(u-w)\}-C_{a_0}q
 +2\Re\langle\partial u,\partial\log f\rangle_\chi.
\]
By Lemma~\ref{lem:algebra}, \(E_{A_u}\ge(\tr A_u)/n\).
Multiplication by \(e^\phi\) proves \eqref{eq:gradient-comparison},
with precisely the dependencies stated in the lemma.
\end{proof}

\begin{lemma}\label{lem:gradient-mixed}
Let \(u,w\) be smooth real functions, with \(u\) strictly admissible
and \(\|u-w\|_\infty\le\beta^{-1}\). Use the weights
\eqref{eq:gradient-weight}, with \(a_0,\beta>0\) and \(K_{a_0}\ge1\).
In a unitary frame diagonalizing \(A_u\), write \(E_{A_u}=\sum_i F^{i \bar i}\lambda_i^2\).
Then
\begin{equation}\label{eq:gradient-lambda-mixed}
\begin{aligned}
 2\left|\Re\sum_i F^{i \bar i}\lambda_i u_i\phi_{\bar i}\right|
 &\le \tfrac18E_{A_u}+
 C_{a_0}\bigl[|\partial u|_\chi^2+(1+\|\partial w\|_\infty)^2\bigr]\Gamma(u-w)\\
 &\quad+C_{a_0}(1+\|\partial w\|_\infty)^4\mathcal F.
\end{aligned}
\end{equation}
\begin{equation}\label{eq:gradient-background-mixed}
\begin{aligned}
 2\left|\Re\sum_{i,p}F^{i \bar i}\omega_{i\bar p}u_p\phi_{\bar i}\right|
 &\le\tfrac14q\Gamma(u-w)\\
 &\quad+C_{a_0}\bigl[1+(1+\|\partial w\|_\infty)|\partial u|_\chi\bigr]\mathcal F.
\end{aligned}
\end{equation}
Here \(C_{a_0}\) depends only on \(a_0,n,\|\omega\|_{L^\infty(\chi)}\);
it is independent of \(\beta,K_{a_0}\) and derivatives of \(u,w\).
\end{lemma}
\begin{proof}
The comparison radius gives \(|-a_0+2\beta(u-w-\beta^{-1})|\le a_0+4\).
Hence
\[
\begin{gathered}
 \phi_i=[-a_0+2\beta(u-w-\beta^{-1})](u-w)_i-a_0w_i,\\
 \Gamma(u)\le2\Gamma(u-w)+2(1+\|\partial w\|_\infty)^2\mathcal F.
\end{gathered}
\]
The last inequality follows from \(u_i=(u-w)_i+w_i\) and
\(\Gamma(w)\le (1+\|\partial w\|_\infty)^2\mathcal F\).

For \eqref{eq:gradient-lambda-mixed}, weighted Cauchy--Schwarz
and \(2xy\le x^2/8+8y^2\) give
\[
 2\left|\Re\sum_i F^{i \bar i}\lambda_i u_i\phi_{\bar i}\right|
 \le \tfrac18E_{A_u}+8\sum_i F^{i \bar i}|u_i|^2|\phi_i|^2.
\]
Using the formula for \(\phi_i\), we estimate
\[
\begin{aligned}
 \sum_i F^{i \bar i}|u_i|^2|\phi_i|^2
 &\le2(a_0+4)^2|\partial u|_\chi^2\Gamma(u-w)
       +2a_0^2(1+\|\partial w\|_\infty)^2\Gamma(u)\\
 &\le\bigl[2(a_0+4)^2|\partial u|_\chi^2+4a_0^2(1+\|\partial w\|_\infty)^2\bigr]\Gamma(u-w)\\
 &\quad+4a_0^2(1+\|\partial w\|_\infty)^4\mathcal F.
\end{aligned}
\]
This proves the first estimate.

For \eqref{eq:gradient-background-mixed}, set
\(b_i=\sum_p\omega_{i\bar p}u_p\) and
\(C_\omega=\sup_M|\chi^{-1}\omega|_{\mathrm{op}}\). Then
\[
 \sum_i F^{i \bar i}|b_i|^2
 \le C_\omega^2|\partial u|_\chi^2\mathcal F.
\]
Split \(\phi_{\bar i}=[-a_0+2\beta(u-w-\beta^{-1})] (u-w)_{\bar i}-a_0w_{\bar i}\).
Weighted Cauchy--Schwarz and Young's inequality yield
\[
\begin{aligned}
 &2\left|\Re\sum_i F^{i \bar i}[-a_0+2\beta(u-w-\beta^{-1})]b_i(u-w)_{\bar i}\right|\\
 &\quad\le2(a_0+4)
       \left(\sum_i F^{i \bar i}|b_i|^2\right)^{1/2}\Gamma(u-w)^{1/2}\\
 &\quad\le\tfrac14q\Gamma(u-w)
       +\frac{4(a_0+4)^2}{q}\sum_i F^{i \bar i}|b_i|^2\\
 &\quad\le\tfrac14q\Gamma(u-w)
       +4(a_0+4)^2C_\omega^2\mathcal F,
\end{aligned}
\]
where the last line uses \(|\partial u|_\chi^2/q\le1\).
The remaining term satisfies
\[
\begin{aligned}
 2a_0\left|\Re\sum_i F^{i \bar i} b_iw_{\bar i}\right|
 &\le2a_0
       \left(\sum_i F^{i \bar i}|b_i|^2\right)^{1/2}\Gamma(w)^{1/2}\\
 &\le2a_0C_\omega (1+\|\partial w\|_\infty)|\partial u|_\chi\mathcal F.
\end{aligned}
\]
Adding these bounds proves \eqref{eq:gradient-background-mixed}.
\end{proof}

\begin{lemma}\label{lem:gradient-energy}
Under the hypotheses of Lemma~\ref{lem:gradient-bochner}, suppose
\(\log f\in C^\alpha(M)\), where \(f=\sigma_k(A_u)\). For
\(d\nu=f\chi^n\), set
\begin{align*}
 I_p&=\int W^p\,d\nu,&
 E_p&=\int W^{p-2}\Gamma(W)\,d\nu,\\
 J_p&=\int W^p q^{-1}\tr A_u\,d\nu.
\end{align*}
For every \(p\ge2\),
\begin{equation}\label{eq:gradient-energy}
 \begin{aligned}
 (p-1)E_p+J_p+\int W^p\{\mathcal F+\Gamma(u-w)\}\,d\nu
 &\le Cp^{(1-\alpha)/(2\alpha)} I_p.
 \end{aligned}
\end{equation}
Here \(C\) depends only on \(M,\chi,\omega,n,k,C_0,f_-,f_+,\alpha\)
and \(\|\log f\|_{C^\alpha(M)}\); it is independent of \(p\)
and derivatives of the comparison function \(w\).
\end{lemma}
\begin{proof}
Multiply \eqref{eq:gradient-comparison} by \(W^{p-1}\) and integrate.
The symmetric identity \eqref{eq:ibp} gives
\[
 \int_M W^{p-1}LW\,d\nu
 =-(p-1)\int_M W^{p-2}\Gamma(W)\,d\nu=-(p-1)E_p.
\]
Hence
\begin{equation}\label{eq:gradient-integrated}
\begin{aligned}
 &(p-1)E_p+cJ_p
       +c\int_M W^p\{\mathcal F+\Gamma(u-w)\}\,d\nu\\
 &\qquad\le CI_p
       -2\int_M\frac{W^p}{q}
                   \Re\langle\partial u,\partial f\rangle_\chi\,dV.
\end{aligned}
\end{equation}
Here \(dV=\chi^n\). The forcing term has this form because
\(e^\phi W^{p-1}=W^p/q\), \(d\nu=f\,dV\), and
\(f\,\partial\log f=\partial f\).

For \(0<r\le1\), define the heat regularization and its remainder by
\[
 f_r:=e^{r^2\Delta_\chi}f,\qquad b_r:=f-f_r.
\]
Thus \(f_r\) is the smoothing of \(f\) at spatial scale \(r\).
Write \(f=f_r+b_r\), where the standard heat estimates give
\(\|b_r\|_\infty\le Cr^\alpha\) and
\(\|\partial f_r\|_\infty\le Cr^{\alpha-1}\).
Since \(|\partial u|_\chi/q\le1\) and \(dV,d\nu\) are uniformly comparable,
\[
 \left|-2\int_M\frac{W^p}{q}
       \Re\langle\partial u,\partial f_r\rangle_\chi\,dV\right|
 \le Cr^{\alpha-1}I_p.
\]
For the remainder, apply \eqref{eq:fixed-ibp} to \(u\) with test
function \(b_rW^p/q\). The product rule gives
\[
\begin{aligned}
 \int_M b_r\frac{W^p}{q}\Delta_\chi u\,dV
 ={}&-\int_M\frac{W^p}{q}
             \Re\langle\partial b_r,\partial u\rangle_\chi\,dV\\
 &-\int_M b_r\Re\left\langle
             \partial\!\left(\frac{W^p}{q}\right),\partial u
                         \right\rangle_\chi\,dV.
\end{aligned}
\]
Since \(W=e^\phi q\), the weight satisfies
\[
 \frac{W^p}{q}=e^\phi W^{p-1},\qquad
 \partial\!\left(\frac{W^p}{q}\right)
 =\frac{W^p}{q}\partial\phi
       +(p-1)\frac{W^{p-1}}q\partial W.
\]
Rearranging the first identity and substituting this derivative yields
\begin{equation}\label{eq:gradient-remainder-ibp}
\begin{aligned}
 &-2\int_M W^pq^{-1}
                \Re\langle\partial u,\partial b_r\rangle_\chi\,dV\\
 &\quad=2\int_M b_r\bigl\{W^pq^{-1}\Delta_\chi u
       +W^pq^{-1}\Re\langle\partial\phi,\partial u\rangle_\chi\\
 &\hspace{38mm}+(p-1)W^{p-1}q^{-1}
                 \Re\langle\partial W,\partial u\rangle_\chi\bigr\}\,dV.
\end{aligned}
\end{equation}
The first two terms in \eqref{eq:gradient-remainder-ibp} are bounded by
\(C\|b_r\|_\infty(J_p+I_p)\). For the third term,
Lemma~\ref{lem:algebra} gives
\[
 |\langle\partial W,\partial u\rangle_\chi|
 \le\{\Gamma(W)(\tr A_u)|\partial u|_\chi^2\}^{1/2}.
\]
Young's inequality bounds its integral by
\((p-1)E_p/2+C(p-1)\|b_r\|_\infty^2J_p\).
The first term is absorbed into the energy on the left.
Choose \(r=r_*p^{-1/(2\alpha)}\), with fixed small \(r_*\),
to absorb the remaining \(J_p\) terms. The smooth forcing then
contributes \(Cp^{(1-\alpha)/(2\alpha)}I_p\), proving
\eqref{eq:gradient-energy}.
\end{proof}

\begin{lemma}\label{lem:gradient-iteration}
Assume the hypotheses of Lemma~\ref{lem:gradient-energy}, and use its
comparison function and weight. Then \(\|W\|_\infty\le C\). Besides the dependencies of Lemma~\ref{lem:gradient-energy},
\(C\) depends only on \(\|\partial w\|_\infty\).
\end{lemma}
\begin{proof}
On \( \{ |\partial u|_\chi^2\ge\max(K_0,4(1+\|\partial w\|_\infty)^2) \}\),
the condition \(|\partial u|_\chi\ge2(1+\|\partial w\|_\infty)\), together with
\(|\partial w|_\chi\le1+\|\partial w\|_\infty\), gives
\[
 |\partial(u-w)|_\chi
 \ge|\partial u|_\chi-|\partial w|_\chi
 \ge\tfrac12|\partial u|_\chi.
\]
Lemma~\ref{lem:algebra} gives \(F^{i \bar i}\ge(\tr A_u)^{-1}\), and hence
\[
 \Gamma(u-w)=\sum_i F^{i \bar i}|(u-w)_i|^2
 \ge\frac{|\partial(u-w)|_\chi^2}{\tr A_u}
 \ge\frac{|\partial u|_\chi^2}{4\tr A_u}.
\]
The other condition, \(|\partial u|_\chi^2\ge K_0\), implies
\(q=|\partial u|_\chi^2+K_0\le2|\partial u|_\chi^2\).
Moreover, \eqref{eq:dual-growth} gives
\[
 \frac1{\tr A_u}\ge\frac{k^{k-1}}{nf_+}
                         (\mathcal F)^{-(k-1)}.
\]
Here \(f_+\) is the upper bound for the normalized density
\(\Sk(A_u)\). Combining these inequalities yields
\[
 \Gamma(u-w)\ge\frac{|\partial u|_\chi^2}{4\tr A_u}
 \ge\frac{q}{8\tr A_u}
 \ge\frac{k^{k-1}}{8nf_+}\frac{q}{(\mathcal F)^{k-1}}.
\]
For \(t>0\), differentiation shows that \(t+cq\,t^{-(k-1)}\)
is minimized at \(t=((k-1)cq)^{1/k}\). Its minimum is a
positive multiple of \(q^{1/k}\). Apply this with \(t=\mathcal F\):
\[
 \mathcal F+\Gamma(u-w)\ge cq^{1/k}\ge cW^{1/k}.
\]
The last inequality uses \(0\le\phi\le a_0C_0+4/\beta\).
On  \( \{ |\partial u|_\chi^2 < \max(K_0,4(1+\|\partial w\|_\infty)^2) \}\), \(W\) is bounded by \(C\).
Consequently \eqref{eq:gradient-energy} yields
\[
 I_{p+1/k}\le C(1+p^{(1-\alpha)/(2\alpha)})I_p.
\]
H\"older's inequality, applied between exponents \(p\) and
\(p+1/k\), therefore bounds every fixed \(I_p\) for
\(p\ge2\).

Again by \(F^{i \bar i}\ge (\tr A_u)^{-1}\) from Lemma \ref{lem:algebra},
\[
 \int|\partial W^{p-1/2}|_\chi\,d\nu
 \le CpE_p^{1/2}J_p^{1/2}
 \le Cp^{1/(2\alpha)}I_p.
\]
On the fixed compact real \(2n\)-dimensional manifold, the
\(W^{1,1}\) Sobolev inequality and the density bounds give
\[
 \|h\|_{L^{2n/(2n-1)}(\nu)}
 \le C\left(\int|\partial h|_\chi\,d\nu+\int|h|\,d\nu\right).
\]
Take \(h=W^{p-1/2}\). Since \(W\ge1\), its \(L^1\) term is
at most \(I_p\), and the preceding estimate yields
\[
 \left(\int W^{\frac{2n}{2n-1}(p-1/2)}\,d\nu\right)^{(2n-1)/(2n)}
 \le Cp^{1/(2\alpha)}I_p.
\]
Start at a fixed \(p_0>2n+1\) and set
\(p_{j+1}=\frac{2n}{2n-1}(p_j-1/2)\). The logarithmic sums converge,
and \(\prod_jp_j/(p_j-1/2)<\infty\). Iteration gives
\(\|W\|_\infty\le C\).
\end{proof}

\begin{proposition}\label{prop:gradient}
Every smooth normalized solution of \eqref{eq:main} satisfies
\[
 \|u\|_{C^1}\le C(M,\chi,\omega,n,k,\alpha,\|G\|_{C^\alpha(M)}).
\]
\end{proposition}
\begin{proof}
Lemma~\ref{lem:potential-holder} gives a finite \(\beta^{-1}\)-net of
the smooth normalized solution class with a fixed upper bound
for \(\|G\|_{C^\alpha(M)}\) in \(C^0\), with centers
in that class. Their gradients have a finite maximum. Since
\(\beta\) is independent of the derivatives of the center,
Lemma~\ref{lem:gradient-iteration} applies to every solution
with one of these centers. The weight \(e^\phi\) has a fixed
positive lower bound, so $W$ bounds
\(|\partial u|_\chi\) uniformly. Together with the potential
bound, this proves the proposition.
\end{proof}

\begin{remark}
Note that the idea of comparing the solution to another function $w$ to prove $C^1$ estimate for solution to some fully nonlinear equation was developed in \cite{guo2026stability}.
\end{remark}

\begin{corollary}\label{cor:gradient-flow}
Let \(\psi\) be a smooth admissible solution of
\(\psi_t=F(A_\psi)-g\), and let \(w\) be a smooth admissible
comparison function fixed in time. Assume
\[
 -C_0\le\psi\le0,\qquad 0<f_-\le\Sk(A_\psi)\le f_+,\qquad
 \|\psi-w\|_\infty\le\beta^{-1}.
\]
Choose \(a_0,\beta,K_{a_0}\) as in Lemma~\ref{lem:gradient-bochner}
for these bounds, and define \(\phi,q,W\) by
\eqref{eq:gradient-weight} with \(u\) replaced by \(\psi\).
Then
\begin{equation}\label{eq:gradient-flow}
 (\partial_t-L)W\le C(1+|\psi_t|)W+C|\partial g|\sqrt W.
\end{equation}
The constant depends only on \(M,\chi,\omega,n,k,C_0,f_-,f_+\).
\end{corollary}
\begin{proof}
Here \(A_\psi=\chi^{-1}(\omega+\ddc\psi)\), and \(a,L,\Gamma\)
are evaluated along \(\psi\).
The instantaneous elliptic forcing is \(g+\psi_t\), and
\[
 W_t=2e^\phi\Re\langle\partial\psi,\partial\psi_t\rangle_\chi
       +[-a_0+2\beta(\psi-w-\beta^{-1})]\psi_tW.
\]
Subtract \eqref{eq:gradient-comparison}. The
\(\partial\psi_t\) terms cancel, the favorable terms may be discarded,
and \(e^\phi|\partial\psi|\le C\sqrt W\). This proves
\eqref{eq:gradient-flow} on the stated interval.
\end{proof}

\section{$C^{1,\beta}$ estimate along flow}\label{sec:flow}

\begin{definition}\label{def:flows}
For a smooth solution \(u\) of \eqref{eq:main}, the original flow solves
\begin{equation}\label{eq:original-flow}
 \psi_t=F(A_\psi)-P_tG,\qquad \psi(0)=u.
\end{equation}
Its unshifted restart at an original time \(\tau>0\) solves
\begin{equation}\label{eq:restart-flow}
 \partial_s\varphi^\tau=F(A_{\varphi^\tau})-P_sG,
 \qquad \varphi^\tau(0)=\psi(\tau).
\end{equation}
Both equations are considered on their intervals of smooth admissible
existence. 
\end{definition}

\subsection{Differential inequalities}
Consider a smooth admissible path solving \(\psi_t=F(A_\psi)-g\), with
fixed potential, gradient, velocity, and positive density bounds.
Use \(A_\psi(t)=\chi^{-1}(\omega+\ddc\psi(t))\),
so \(\tr A_\psi(t)=\tr_\chi\omega+\Delta_\chi\psi(t)\), and evaluate
\(a,L,\Gamma\) at this endomorphism.
Order the eigenvalues as \(\lambda_1\ge\cdots\ge\lambda_n\), and define
\[
 D=\partial_t-L,\qquad
 H=\sum_{i,p}F^{i \bar i}(|\psi_{pi}|^2+|\psi_{p\bar i}|^2).
\]
Denote $K= 1+  \|\psi\|_\infty+\|\partial\psi\|_\infty+\|\psi_t\|_\infty.$ Define
\[
 \phi(s)=-\tfrac18\log\left(1-\frac{s}{2(1+K^2)}\right),
 \qquad 0\le s\le K^2.
\]
Define the following function which is the same as
\cite[(2.18)]{PhongTo2021}:
\[
 h(x)=D_1e^{-D_2x},\qquad 0\le x\le2K.
\]
The parameters \(D_1,D_2>0\) are fixed below. We have
\[
 h'<0,\qquad h''=D_2(-h'),\qquad
 \frac{(h')^2}{h''}=h\le D_1,\qquad
 \min_{[0,2K]}(-h')=D_1D_2e^{-2KD_2}.
\]

Define
\begin{equation}\label{e def z}
 Z=\lambda_1\exp\{\phi(|\partial\psi|_\chi^2)+h(\psi+K)\}.
\end{equation}

$Z$ will be the function that we use to estimate $\lambda_1$.

\begin{lemma}\label{lem:parabolic-gradient}
Under the standing K\"ahler assumptions, let \(\psi\) be a smooth
strictly admissible solution of \(\psi_t=F(A_\psi)-g\) on
\(M\times[0,T]\).
Then
\begin{equation}\label{eq:parabolic-gradient}
 \begin{aligned}
 D\phi(|\partial\psi|_\chi^2)
 \le{}&-\phi'H-\phi''\Gamma(|\partial\psi|_\chi^2) +C\mathcal F+C|\partial g|_\chi.
 \end{aligned}
\end{equation}
Here \(\phi'\) and \(\phi''\) are evaluated at
\(|\partial\psi|_\chi^2\), and \(C\) depends only on
\[
 n,\quad
 \|\nabla\omega\|_{L^\infty(\chi)},\quad
 \|\mathrm{Rm}(\chi)\|_{L^\infty(\chi)}.
\]
Moreover,
\begin{equation}\label{eq:parabolic-gradient-hessian}
 H\ge\tfrac12\sum_i F^{i \bar i}\lambda_i^2-C_\omega\mathcal F,
 \qquad C_\omega=\|\chi^{-1}\omega\|_{L^\infty(\mathrm{op})}^2.
\end{equation}
\end{lemma}
\begin{proof}
Use the weighted gradient calculation leading to
\cite[(2.37)]{PhongTo2021}, with vanishing torsion and with their
parabolic operator equal to \(-D\). Our \(\phi\) is one quarter of
their logarithmic weight with \(P=1+K^2\). The cited calculation
and the chain rule give
\begin{equation}\label{eq:parabolic-gradient-identity}
\begin{aligned}
 D\phi(|\partial\psi|_\chi^2)
 ={}&\phi'D|\partial\psi|_\chi^2
       -\phi''\Gamma(|\partial\psi|_\chi^2)\\
 ={}&-\phi'H-\phi''\Gamma(|\partial\psi|_\chi^2)\\
 &+2\phi'\Re\sum_{i,p}F^{i \bar i}\psi_{\bar p}\nabla_p\omega_{i\bar i}\\
 &-\phi'\sum_{i,p,q}F^{i \bar i}R_{i\bar i p\bar q}\psi_q\psi_{\bar p}
 -2\phi'\Re\langle\partial\psi,\partial g\rangle_\chi.
\end{aligned}
\end{equation}
The curvature convention is that of Lemma~\ref{lem:kahler-bochner}.
In the cited estimate the forcing derivatives are absorbed into its
constant; here we retain them explicitly. Only spatial derivatives
of \(g\) enter, so its time dependence does not change the calculation.
On \(0\le s\le K^2\),
\[
 \frac1{16(1+K^2)}\le\phi'(s)\le\frac1{8(1+K^2)},
 \qquad \phi''(s)=8\phi'(s)^2.
\]
Consequently, \(\phi'|\partial\psi|_\chi\) and
\(\phi'|\partial\psi|_\chi^2\) are bounded by absolute constants.
Using \(\sum_i F^{i \bar i}=\mathcal F\) in
\eqref{eq:parabolic-gradient-identity} proves
\eqref{eq:parabolic-gradient} with the stated dependencies.
Finally, \(\psi_{i\bar i}=\lambda_i-\omega_{i\bar i}\) and
\(|x-y|^2\ge|x|^2/2-|y|^2\) give
\[
 H\ge\tfrac12\sum_i F^{i \bar i}\lambda_i^2-\sum_i F^{i \bar i}|\omega_{i\bar i}|^2
 \ge\tfrac12\sum_i F^{i \bar i}\lambda_i^2-C_\omega\mathcal F.
\]
This is \eqref{eq:parabolic-gradient-hessian}.
\end{proof}

\begin{lemma}\label{lem:offcritical}
Let \(\psi\) be a smooth admissible solution of
\(\psi_t=F(A_\psi)-g\) on \(M\times[0,T]\). Assume \(K\ge1\) and
\[
 \|\psi\|_\infty+\|\partial\psi\|_\infty+\|\psi_t\|_\infty\le K,
 \qquad 0<f_-\le\Sk(A_\psi)\le f_+.
\]
Then 
\begin{equation}\label{eq:offcritical}
 D\log Z\le-cH-c\lambda_1-\mathcal F
       +C(1+|\partial g|+|\nabla\bar\nabla g|)+4\Gamma(\log Z)
\end{equation}
holds at the place where $\lambda_1$ is a smooth function. Its constants
depend only on \(M,\chi,\omega,n,k,K,f_-,f_+\).
\end{lemma}
\begin{proof}
Under the density hypothesis, Proposition~\ref{prop:newton-maclaurin}
gives
\begin{equation}\label{e den upp low}
 \lambda_1\ge\frac{\tr A_\psi}{n}\ge\Sk(A_\psi)^{1/k}\ge f_-^{1/k}>0,
 \qquad \lambda_1\le\tr A_\psi\le n\lambda_1.
\end{equation}
We use Phong--T\^o's Hou--Ma--Wu calculation
\cite{HMW,PhongTo2021}, specialized to the K\"ahler case and
\(F=\log S_k\), retaining derivatives of the test function.

Around a point where $\lambda_1$ is locally smooth, we apply
\cite[(2.39)]{PhongTo2021}, with their operator equal to \(-D\),
vanishing torsion, and the weighted gradient bound
\eqref{eq:parabolic-gradient}. Keeping the terms involving $g$ explicit
before they are absorbed into the constant in \cite[(2.39)]{PhongTo2021} gives
\begin{equation}\label{eq:phong-to-combined}
\begin{aligned}
- D \log Z
\geq\;&
- F^{k\bar{k}}
\frac{\left|\widetilde{\lambda}_{1,\bar{k}}\right|^2}{\lambda_1^2}
-\frac{1}{\lambda_1}
F^{l\bar{k},\,s\bar{r}}
\nabla_{\bar{1}} (\omega_{\psi})_{\bar{k}l}
\nabla_{1} (\omega_{\psi})_{\bar{r}s}
-C_4\frac{1}{\lambda_1}
F^{k\bar{k}}
\left|\lambda_{1,\bar{k}}\right|
-C_9\mathcal{F}
\\
&\quad
+\frac{1}{16 (1+K^2)}
F^{q\bar{q}}
\left(
\left|\nabla_q\nabla u\right|^2
+
\left|\nabla_q\overline{\nabla}u\right|^2
\right)
+\phi'' F^{q\bar{q}}
\nabla_q|\nabla u|^2
\nabla_{\bar{q}}|\nabla u|^2
\\
&\quad
+h'(\widetilde{v})
\left(
F^{k\bar{k}}
\nabla_k\nabla_{\bar{k}}\widetilde{v}
-\partial_t\widetilde{v}
\right)
+h''(\widetilde{v})
F^{k\bar{k}}
\nabla_k\widetilde{v}
\nabla_{\bar{k}}\widetilde{v} -C \frac{|g_{1 \bar 1}|}{\lambda_1}-C |\partial g|.
\end{aligned}
\end{equation}
Next, we consider the two cases as in \cite{PhongTo2021}:

\textbf{Case 1}: $\theta \lambda_1 \le - \lambda_n$. Phong-To did the calculation at a maximum point of $\log Z$, while we are not. As a result, we can't assume that $\nabla \log Z =0$. Instead,  compared with \cite[(2.45)]{PhongTo2021}, the estimate we get is:
\begin{equation}\label{e Dlogz case1}
-D \log Z \ge \frac{C_1}{K^2}\mathcal F \lambda_1^2 + \frac{C_5}{K^2} \sum_p F^{i \bar i} |\psi_{i p}|^2 - C_2 K^2 \mathcal F - C_3 \Gamma(\log Z)- C_4(1+|\partial g|+ |\partial \bar \partial g|).
\end{equation}
Without loss of generality, we can assume that $\lambda_1$ is big and satisfies $\lambda_1> \sqrt{\frac{2 C_2 K^4}{C_1}}$. In fact, otherwise $|\lambda|$ is bounded, and $\mathcal{F}$ is bounded as a result. Then, (\ref{eq:offcritical}) would follows from (\ref{e Dlogz case1}) immediately.

With (\ref{e den upp low}), we can use \cite[Lemma 9]{Szekelyhidi2018} to get that $\mathcal F \ge \tau$ for some constant $\tau>0$. Combining this with the above lower bound on $\lambda_1$, there exist positive constants $C_4, C_5$ and $C_6$ such that
\begin{equation*}
\mathcal F \lambda_1^2 \ge C_4 \lambda_1 + C_5 \mathcal F - C_6.
\end{equation*}
 
We calculate that:
\begin{equation}\label{e flambda kf}
\begin{split}
\frac{C_1}{K^2}\mathcal F \lambda_1^2 -C_2 K^2 \mathcal F &= \frac{C_1}{4K^2}\mathcal F \lambda_1^2 + \frac{C_1}{4K^2}\mathcal F \lambda_1^2 + \frac{C_1}{2K^2}\mathcal F \lambda_1^2  - C_2 K^2 \mathcal F \\
& \ge \frac{C_1}{4K^2}\mathcal F \lambda_1^2  +  \frac{C_1}{4K^2}(C_4 \lambda_1 + C_5 \mathcal F - C_6) +  C_2 K^2 \mathcal F - C_2 K^2 \mathcal F \\
& =  \frac{C_1}{4K^2}\mathcal F \lambda_1^2  +  \frac{C_1}{4K^2}(C_4 \lambda_1 + C_5 \mathcal F - C_6)
\end{split}
\end{equation}

Combining (\ref{e Dlogz case1}) with (\ref{e flambda kf}), we can conclude the proof of (\ref{eq:offcritical}) in the Case 1.

\textbf{Case 2:} $\theta \lambda_1 > -\lambda_n$.

We add the term $C H$ which is discarded in \cite[(2.63)]{PhongTo2021} to get:
\begin{equation}
\begin{split}
-D \log Z & \ge C H + \frac{1}{2} h'' F^{k \bar k} |\nabla_{\bar k}\psi|^2 + h' (F^{k \bar k} \nabla_k \nabla_{\bar k} \psi- \partial_t \psi) -2 (h')^2 \theta^{-1} F^{1\bar 1} K^2 -C_{\epsilon}|h'| F^{k \bar k} |\nabla_{\bar k}\psi|^2 \\
&-\epsilon C|h'|\mathcal F - C \mathcal F - C(1+|\partial g|+ |\partial \bar \partial g|) - C \Gamma(\log Z),
\end{split}
\end{equation}
where $\theta$ is some constant.

As in \cite{PhongTo2021}, we consider two subcases:

\textbf{Case 2.1:} $F^{k \bar k}(\underline{\psi}_{k \bar k} - \psi_{k \bar k})+(\partial_t \psi- \partial_t \underline{\psi}) \ge \kappa \mathcal F.$ Note that for complex $k$-Hessian equation, we can take $\underline{u}=0$. As in \cite[(2.66)]{PhongTo2021}, in this case, we can get that:
\begin{equation}\label{e -dlogz case2.1}
-D \log Z \ge C H + F^{1 \bar 1} (\frac{\lambda_1^2}{20K^2}- 2 (h')^2 \theta^{-1}K^2)- C_{23} \mathcal F -\frac{1}{2} h' \kappa \mathcal F - C(1+|\partial g|+ |\partial \bar \partial g|)- C \Gamma(\log Z).
\end{equation}

We can choose $D_1$ big such that $-C_{23}- \frac{1}{2} h' \kappa \ge 1$. We can also assume that $\lambda_1$ is big, satisfying $\frac{\lambda_1^2}{40K}- 2 (h')^2 \theta^{-1}K >0$ using the argument in Case 1. Using \cite[Lemma 3.1]{CW}, we have that $\lambda_1 \sigma_{k-1}(\lambda | 1) \ge \frac{k}{n}\sigma_k(\lambda).$ This implies that 
\begin{equation*}
F^{1\bar 1}\lambda_1 = \frac{\lambda_1 \sigma_{k-1}(\lambda| 1)}{\sigma_k(\lambda)}\ge \frac{k}{n}.
\end{equation*}

Using the above calculation, we can derive from (\ref{e -dlogz case2.1}) that:
\begin{equation*}
-D \log Z \ge CH + F^{1\bar 1}\frac{\lambda_1^2}{40K^2} + \mathcal F  - C(1+|\partial g|+ |\partial \bar \partial g|)\ge  CH + \frac{k \lambda_1}{40K^2 n} + \mathcal F  - C(1+|\partial g|+ |\partial \bar \partial g|) - C \Gamma(\log Z)
\end{equation*}

\textbf{Case 2.2} $F^{1 \bar 1}\ge \kappa \mathcal F$.
In this case, we can follow \cite[(2.68)]{PhongTo2021} to get:
\begin{equation*}
-D \log Z  \ge C H + \kappa \mathcal F (\frac{\lambda_1^2}{20K^2} -2(h')^2 \theta^{-1}K^2) -C_{24}\mathcal F - h' F^{k \bar k} (\omega_{\psi})_{k \bar k} - C(1+|\partial g|+ |\partial \bar \partial g|) - C \Gamma(\log Z)
\end{equation*}

Note that $- h' F^{k \bar k} (\chi_\psi)_{k \bar k}\ge 0 $. We can throw this term away. Then, we can use the argument in the Case 1 to use the term $\mathcal F \lambda_1^2$ to generate all the good terms that we need. This concludes the proof of this Proposition.
\end{proof}

\subsection{$C^0$ and $C^1$ estimate along flow}

\begin{lemma}[Velocity and density bounds]\label{lem:flow-density}
Let \(\psi\) and \(\varphi^\tau\) be the original flow and its
unshifted restarts from Definition~\ref{def:flows}, starting from a
smooth solution \(u\) of \eqref{eq:main}. For any \(0\le t,s\le1\) such that $\psi(t)$ and $\varphi^{\tau}$ exist,
\begin{equation}\label{eq:flow-velocity-bounds}
 \|\psi_t(t)\|_\infty\le Ct^{\alpha/2},\qquad
 \|\partial_s\varphi^\tau(s)\|_\infty
       \le C(\tau^{\alpha/2}+s^{\alpha/2}).
\end{equation}
There are constants \(0<f_-\le f_+<\infty\) such that
\begin{equation}\label{eq:flow-density-window}
 \begin{aligned}
 f_-&\le\Sk(A_\psi(t))\le f_+,\\
 f_-&\le\Sk(A_{\varphi^\tau}(s))\le f_+.
 \end{aligned}
\end{equation}
The constants \(C,f_-,f_+\) depend only on
\(M,\chi,n,\alpha,\|G\|_{C^\alpha(M)}\).
\end{lemma}
\begin{proof}
Differentiate the original flow in time, with \(L\) evaluated along
that flow. Since \(u\) solves \eqref{eq:main},
\[
 (\partial_t-L)\psi_t=-\partial_tP_tG,\qquad
 \psi_t(0)=\log\Sk(A_u)-G=0.
\]
The maximum principle and Lemma~\ref{lem:heat-regularization} give
\[
 \|\psi_t(t)\|_\infty
 \le\int_0^t\|\partial_rP_rG\|_\infty\,dr
 \le C\int_0^t r^{\alpha/2-1}\,dr
 \le Ct^{\alpha/2}.
\]
At the start of a restart,
\[
 \partial_s\varphi^\tau(0)
 =\log\Sk(A_\psi(\tau))-G
 =\psi_t(\tau)+P_\tau G-G.
\]
Its supremum norm is at most \(C\tau^{\alpha/2}\) by the preceding
bound and \eqref{eq:heat-envelope}. Along the restart,
\[
 (\partial_s-L)\partial_s\varphi^\tau=-\partial_sP_sG,
\]
where \(L\) is now evaluated along \(\varphi^\tau\). The same
maximum-principle argument proves the second bound in
\eqref{eq:flow-velocity-bounds}.

The flow equations imply
\[
 \Sk(A_\psi(t))=\exp\bigl(P_tG+\psi_t(t)\bigr),\qquad
 \Sk(A_{\varphi^\tau}(s))
   =\exp\bigl(P_sG+\partial_s\varphi^\tau(s)\bigr).
\]
Since \(\|P_tG\|_\infty\le\|G\|_\infty\) and \(t,s,\tau\le1\),
we may take
\[
 f_-=\exp(-\|G\|_\infty-2C),\qquad
 f_+=\exp(\|G\|_\infty+2C).
\]
This proves \eqref{eq:flow-density-window}.
\end{proof}

\begin{lemma}[Gradient propagation from a comparison center]\label{lem:comparison-gradient-propagation}
Let \(\psi\) be a smooth admissible solution of
\(\psi_t=F(A_\psi)-g\) on \(M\times[0,T]\), with \(T\le1\). Assume
\[
 -C_0\le\psi\le0,\qquad
 0<f_-\le\Sk(A_\psi)\le f_+.
\]
Choose \(a_0,\beta,K_{a_0}\) as in Lemma~\ref{lem:gradient-bochner}
for these potential and density bounds.
Let \(w\) be a fixed smooth admissible function such that
\(\|\psi(t)-w\|_\infty\le\beta^{-1}\) for \(0\le t\le T\),
and define \(\phi,q,W\) as in \eqref{eq:gradient-weight}.
Then
\begin{equation}\label{eq:comparison-gradient-propagation}
 \|\sqrt{W(t)}\|_\infty
 \le
 \exp\!\left(C\int_0^t(1+\|\psi_t(s)\|_\infty)\,ds\right)
 \left(\|\sqrt{W(0)}\|_\infty
       +C\int_0^t\|\partial g(s)\|_\infty\,ds\right).
\end{equation}
The constant depends only on \(M,\chi,\omega,n,k,C_0,f_-,f_+\).
\end{lemma}
\begin{proof}
Corollary~\ref{cor:gradient-flow} gives
\[
 (\partial_t-L)W\le C(1+|\psi_t|)W+C|\partial g|\sqrt W.
\]
At a spatial maximum of \(W\), one has \(LW\le0\).
Consequently, for \(m(t)=\max_M W(t)>0\), the upper right Dini
derivative satisfies
\[
 D^+m(t)\le C(1+\|\psi_t(t)\|_\infty)m(t)
                  +C\|\partial g(t)\|_\infty\sqrt{m(t)}.
\]
Dividing by \(2\sqrt{m(t)}\) and applying Gronwall's inequality
proves \eqref{eq:comparison-gradient-propagation}.
\end{proof}

We record $C^0$ and $C^1$ estimate for $\psi$ and $\varphi^{\tau}$:
\begin{proposition}[Uniform flow bounds and comparison]\label{prop:flow-bounds}
Let \(u\) be a smooth normalized solution of \eqref{eq:main},
and let \(\psi,\varphi^\tau\) be the flows of Definition~\ref{def:flows}.
There are \(T_*\le1\), \(0<\tau_0\le T_*/2\) and constants $C, f_+, f_-$, depending only on
\[
 M,\chi,\omega,n,k,\alpha,\|G\|_{C^\alpha(M)},
\]
such that for any \(0\le t,s\le T_*\) and \(0<\tau\le\tau_0\) which make $\psi(t)$ and $\varphi^{\tau}(s)$ exist, we have the a priori estimate:
\begin{equation}\label{eq:common-flow-uniform-bounds}
\begin{gathered}
 \|\psi(t)\|_\infty+\|\partial\psi(t)\|_\infty
       +\|\psi_t(t)\|_\infty\le C,\\
 \|\varphi^\tau(s)\|_\infty+\|\partial\varphi^\tau(s)\|_\infty
       +\|\partial_s\varphi^\tau(s)\|_\infty\le C,\\
 0<f_-\le\Sk(A_\psi(t))\le f_+,\qquad
 0<f_-\le\Sk(A_{\varphi^\tau}(s))\le f_+.
\end{gathered}
\end{equation}
Moreover,
\begin{gather}
 \|\psi_t(t)\|_\infty\le Ct^{\alpha/2},\qquad
 \|\psi(t)-u\|_\infty\le Ct^{1+\alpha/2},\label{eq:displacement}\\
 \|\varphi^\tau(s)-\psi(s)\|_\infty\le C\tau^{1+\alpha/2},
 \qquad 0\le s\le T_*,\quad0<\tau\le\tau_0.\label{eq:restart-comparison}
\end{gather}
\end{proposition}
\begin{proof}
Lemma~\ref{lem:flow-density} gives \eqref{eq:flow-velocity-bounds}
and \eqref{eq:flow-density-window} on the respective existence
intervals. Integrating the original velocity estimate gives
\[
 \|\psi(t)-u\|_\infty
 \le C\int_0^t r^{\alpha/2}\,dr
 \le Ct^{1+\alpha/2},
\]
which proves \eqref{eq:displacement}. Along a restart,
\[
\begin{aligned}
 \|\varphi^\tau(s)-u\|_\infty
 &\le\|\psi(\tau)-u\|_\infty
       +\int_0^s\|\partial_r\varphi^\tau(r)\|_\infty\,dr\\
 &\le C\bigl(\tau^{1+\alpha/2}
             +\tau^{\alpha/2}s+s^{1+\alpha/2}\bigr).
\end{aligned}
\]
Together with the uniform potential bound for \(u\), these
estimates bound both potentials uniformly for \(t,s,\tau\le1\).
Choose the comparison parameter \(\beta\) in
Lemma~\ref{lem:comparison-gradient-propagation} using these
potential bounds and the fixed density window. We may then fix
\(T_*\le1\) and \(0<\tau_0\le T_*/2\) so that
\[
 \|\psi(t)-u\|_\infty\le\beta^{-1},\qquad
 \|\varphi^\tau(s)-u\|_\infty\le\beta^{-1}
\]
throughout the stated ranges, whenever the paths exist.

Apply Lemma~\ref{lem:comparison-gradient-propagation} first to
the original flow with comparison center \(u\). Its initial
gradient is bounded by Proposition~\ref{prop:gradient}, and
Lemma~\ref{lem:heat-regularization} gives
\[
 \int_0^t\|\partial P_rG\|_\infty\,dr
 \le Ct^{(1+\alpha)/2}.
\]
The velocity integral is also uniformly bounded by
\eqref{eq:flow-velocity-bounds}. Thus
\(\|\partial\psi(t)\|_\infty\le C\).
In particular, every defined restart begins with
\(\|\partial\varphi^\tau(0)\|_\infty
=\|\partial\psi(\tau)\|_\infty\le C\).
Applying the same lemma to each restart with comparison center
\(u\) now gives \(\|\partial\varphi^\tau(s)\|_\infty\le C\).
The potential, velocity, and density bounds already obtained
complete \eqref{eq:common-flow-uniform-bounds}.

Finally, the original flow and the restart have the same forcing
\(P_sG\). The scalar comparison principle and
\eqref{eq:displacement} imply
\[
 \|\varphi^\tau(s)-\psi(s)\|_\infty
 \le\|\psi(\tau)-u\|_\infty
 \le C\tau^{1+\alpha/2}
\]
on their common existence interval. This proves
\eqref{eq:restart-comparison}.
\end{proof}

\subsection{Common existence intervals}

\begin{lemma}[Multiplicative trace propagation]\label{lem:coarse-propagation}
Let \(\psi\) be a smooth admissible solution of
\(\psi_t=F(A_\psi)-P_tG\) on \(M\times[0,T]\), with \(0<T\le1\)
and \(G\in C^\alpha(M)\). Assume the potential, gradient, velocity,
and density bounds of Lemma~\ref{lem:offcritical}, with constants
\(K,f_-,f_+\). Then,
\begin{equation}\label{eq:coarse-propagation}
 \sup_M\tr A_\psi(t)\le C\sup_M\tr A_\psi(0).
\end{equation}
The constant depends only on
\(M,\chi,\omega,n,k,K,f_-,f_+,\alpha,\|G\|_{C^\alpha(M)}\).
The same estimates hold for each unshifted restart in its own time
variable whenever these hypotheses hold.
\end{lemma}
\begin{proof}
At a spatial maximum of \(Z\), one has
\(\Gamma(\log Z)=0\) and \(L\log Z\le0\).
Discarding the negative terms in \eqref{eq:offcritical} therefore gives
\begin{equation}\label{e D log z}
 D^+\!\left(\max_M\log Z(t)\right)
 \le C\bigl(1+\|\partial P_tG\|_\infty
                +\|\nabla\bar\nabla P_tG\|_\infty\bigr),
 \qquad t>0,
\end{equation}
where \(D^+\) is the upper right Dini derivative. If the maximum eigenvalue is not smooth at the maximum point above, then we can take a normal coordinate of $\chi$ around the maximum point such that $A_{\psi}$ is diagonalized at that point and $(A_{\psi})_{11}=\lambda_1$. Then, we take $B=diag (0,\epsilon,\ldots, \epsilon )$. We can replace $A_{\psi}$ by $A_{\psi}-B$. Denote $\widetilde{\lambda}_1 \ge \ldots \ge \widetilde{\lambda}_n$ as the eigenvalues of $A_{\psi}-B$. Then, $\widetilde{\lambda}_1$ is smooth and at the maximum point, $\widetilde{\lambda}_1 = \lambda_1$. The computation will be the same as in the case where $\lambda_1$ is smooth. As a result, we can just assume that $\lambda_1$ is smooth at the maximum point.

Integrating (\ref{e D log z}) yields
\[
 \sup_M Z(t)\le\sup_M Z(0)\,
 \exp\!\left(C\int_0^t
   [1+\|\partial P_sG\|_\infty
       +\|\nabla\bar\nabla P_sG\|_\infty]\,ds\right).
\]
By Lemma~\ref{lem:heat-regularization},
\[
 \int_0^t
   [1+\|\partial P_sG\|_\infty
       +\|\nabla\bar\nabla P_sG\|_\infty]\,ds
 \le C\bigl(t+t^{(1+\alpha)/2}+t^{\alpha/2}\bigr)
 \le Ct^{\alpha/2}.
\]
The integral is finite because \(\alpha>0\); continuity gives the
endpoint at \(t=0\). 
Using \(Z\asymp\tr A_\psi\), with the comparison constants of
Lemma~\ref{lem:offcritical}, gives \eqref{eq:coarse-propagation}.
An unshifted restart has the same forcing \(P_sG\), so the argument
applies in its own clock.
\end{proof}

\begin{proposition}[Common existence intervals]\label{prop:common-flow}
For every smooth normalized solution \(u\) of \eqref{eq:main},
the original flow and every unshifted restart with
\(0<\tau\le\tau_0\) exist on the common closed interval
\([0,T_*]\), where \(T_*,\tau_0\) are those chosen in
Proposition~\ref{prop:flow-bounds}. In particular, these times
depend only on
\[
 M,\chi,\omega,n,k,\alpha,\|G\|_{C^\alpha(M)}.
\]
\end{proposition}
\begin{proof}
For each fixed smooth \(u,G\), strict admissibility and standard
parabolic local existence give an original solution on a short
interval. Proposition~\ref{prop:flow-bounds} bounds its potential,
gradient, velocity, and density as long as it exists with
\(t\le T_*\). Lemma~\ref{lem:coarse-propagation} therefore gives
\[
 \sup_M\tr A_\psi(t)\le C\sup_M\tr A_u.
\]
For this continuation argument the right-hand side may depend
on the individual \(C^2\) norm of \(u\). Lemma~\ref{lem:algebra}
bounds the complex Hessian, and the positive density floor and
the Newton--Maclaurin inequalities in
Proposition~\ref{prop:newton-maclaurin} keep its eigenvalues in
a compact subset of \(\Gamma_k\).

At any putative maximal time not exceeding \(T_*\), these
bounds, together with the technique that view a partial differential equation which uses complex Hessian of the solution as a partial differential equation using real Hessian of the solution (see \cite{TWWY2015}),  verify the compact-range hypotheses of
\cite[Theorem~5.1]{Chu}. On a time interval bounded away from
zero, that estimate and Schauder bootstrap give full real
parabolic regularity and a smooth strictly admissible limit
 for the original flow and the unshifted restart with $0< \tau \le \tau_0$ on their existence intervals intersected with $[0, T_*]$. This a priori estimate ensures that the original flow and every unshifted restart with $0< \tau \le \tau_0$ exist on $[0, T_*]$.
\end{proof}

\subsection{Entropy and trace decay}

The following estimates are written for the original potential \(\psi\).
For a restart, replace
\(\psi(t)\), \(A_\psi(t)\), and \(\psi_t(t)\) by
\(\varphi^\tau(s)\), \(A_{\varphi^\tau}(s)\), and
\(\partial_s\varphi^\tau(s)\), respectively.
The operators and measure are evaluated on that same path.

\begin{lemma}[Weighted energy for the evolving measure]\label{lem:weighted-energy}
On a path of Proposition~\ref{prop:common-flow}, suppose \(|\psi_t|\le B\)
on the spacetime interval. Let \(w>0\) be smooth and satisfy
\(D\log w\le q+\gamma\Gamma(\log w)\).
For \(Y_p=\int e^{\psi_t^2}w^p\,d\nu\), where
\(d\nu=\Sk(A_\psi)\chi^n\), one has
\begin{equation}\label{eq:weighted-energy}
 \begin{aligned}
 Y_p'+c\int e^{\psi_t^2}w^p
       [p^2\Gamma(\log w)+\Gamma(\psi_t)]\,d\nu
 &\le p\int q\,e^{\psi_t^2}w^p\,d\nu
       +C\|g_t\|_\infty Y_p,\qquad p\ge p_0 .
 \end{aligned}
\end{equation}
The constants \(c,C>0\) depend only on \(B\). For a given \(\gamma\)
one may take \(p_0=\max\{2,\gamma(4B^2+2B+2)\}\).
\end{lemma}
\begin{proof}
Put \(W=e^{\psi_t^2}w^p\), so \(Y_p=\int_MW\,d\nu\).
Since \(d\nu=e^{F(A_\psi)}\chi^n\), \(F(A_\psi)=\psi_t+g\),
and \(D\psi_t=-g_t\), the measure satisfies
\[
 \partial_t(d\nu)
 =(\psi_{tt}+g_t)d\nu=(L\psi_t)d\nu.
\]
Differentiating the integral and using \eqref{eq:ibp}, we obtain
\begin{equation}\label{eq:evolving-weighted-integral}
\begin{aligned}
 Y_p'
 &=\int_M\bigl(DW+LW+W L\psi_t\bigr)\,d\nu\\
 &=\int_M DW\,d\nu-\int_M\Gamma(W,\psi_t)\,d\nu.
\end{aligned}
\end{equation}
Here \(\int_M LW\,d\nu=0\), while integration by parts converts
\(\int_MW L\psi_t\,d\nu\) into the second term.

The chain rule gives
\[
 DW=W\bigl(D\log W-\Gamma(\log W)\bigr).
\]
Using \(\log W=\psi_t^2+p\log w\), we calculate
\[
\begin{aligned}
 D\log W
 &=2\psi_tD\psi_t-2\Gamma(\psi_t)+pD\log w\\
 &\le-2\psi_tg_t-2\Gamma(\psi_t)
       +pq+p\gamma\Gamma(\log w),\\
 \Gamma(\log W)
 &=4\psi_t^2\Gamma(\psi_t)+p^2\Gamma(\log w)
       +4p\psi_t\Gamma(\log w,\psi_t).
\end{aligned}
\]
Subtracting the second expression from the first yields
\[
\begin{aligned}
 DW\le W\bigl[&
 pq-p(p-\gamma)\Gamma(\log w)
       -4p\psi_t\Gamma(\log w,\psi_t)\\
 &-(4\psi_t^2+2)\Gamma(\psi_t)-2\psi_tg_t\bigr].
\end{aligned}
\]
Moreover, \(dW=W\,d\log W\) implies
\[
 \Gamma(W,\psi_t)
 =W\bigl[2\psi_t\Gamma(\psi_t)
             +p\Gamma(\log w,\psi_t)\bigr].
\]
Substituting these two formulas into
\eqref{eq:evolving-weighted-integral} gives
\begin{equation}\label{eq:weighted-energy-expansion}
\begin{aligned}
 Y_p'\le\int_M W\bigl[&
 pq-p(p-\gamma)\Gamma(\log w)\\
 &-p(4\psi_t+1)\Gamma(\log w,\psi_t)\\
 &-(4\psi_t^2+2\psi_t+2)\Gamma(\psi_t)
       -2\psi_tg_t\bigr]\,d\nu.
\end{aligned}
\end{equation}
In the covectors \(p\,\partial\log w\) and \(\partial\psi_t\), paired
using \(a^{i\bar j}\), the quadratic form subtracted in
\eqref{eq:weighted-energy-expansion} has coefficient matrix
\[
 \begin{pmatrix}
  1-\gamma/p & (4\psi_t+1)/2\\
  (4\psi_t+1)/2 & 4\psi_t^2+2\psi_t+2
 \end{pmatrix}.
\]
Its determinant is
\[
 \left(1-\frac{\gamma}{p}\right)(4\psi_t^2+2\psi_t+2)
    -\frac{(4\psi_t+1)^2}{4}
 =\frac74-\frac{\gamma}{p}(4\psi_t^2+2\psi_t+2).
\]
Since \(|\psi_t|\le B\) and \(p\ge p_0\), the determinant is at least
\(3/4\). The lower right entry is positive, so the matrix is positive
definite. Its trace is at most \(4B^2+2B+3\), hence its smaller
eigenvalue is at least \(3/[4(4B^2+2B+3)]\). The quadratic form thus
bounds \(c[p^2\Gamma(\log w)+\Gamma(\psi_t)]\) from below.
Using \(|2\psi_t g_t|\le 2B\|g_t\|_\infty\) in
\eqref{eq:weighted-energy-expansion} proves \eqref{eq:weighted-energy}.
\end{proof}

\begin{lemma}[Entropy dissipation]\label{lem:entropy}
On the original flow and each restart of Proposition~\ref{prop:common-flow},
there is a fixed \(p_0>4\) such that, for
\(Y_p=\int e^{\psi_t^2}Z^p\,d\nu\),
\begin{align}
 Y_p'&+cp^2\int e^{\psi_t^2}Z^p\Gamma(\log Z)\,d\nu
       +c\int e^{\psi_t^2}Z^p\Gamma(\psi_t)\,d\nu\notag\\
 &\quad+cp\int e^{\psi_t^2}Z^p(H+\lambda_1)\,d\nu
       \le Cp\rho(t)Y_p,\qquad p\ge p_0,\label{eq:largest-dissipation}
\end{align}
where \(d\nu=e^{F(A_\psi)}\chi^n\) and
\(\rho(t)=1+\|\partial g\|_\infty+\|\nabla\bar\nabla g\|_\infty
+\|g_t\|_\infty\le C(1+t^{\alpha/2-1})\).
Moreover, the dissipation controls both gradient energies:
\begin{equation}\label{eq:eigenvalue-energy-control}
\begin{aligned}
 &\int_M\left[
 \Gamma(Z^{p/2})+
 \Gamma\!\left(e^{\psi_t^2/2}Z^{p/2}\right)\right]\,d\nu\\
 &\qquad\le C\int_M e^{\psi_t^2}Z^p
 \left[p^2\Gamma(\log Z)+\Gamma(\psi_t)\right]\,d\nu,
 \qquad p\ge p_0.
\end{aligned}
\end{equation}
The constants \(p_0,c,C\) depend only on
\(M,\chi,\omega,n,k,\alpha,\|G\|_{C^\alpha(M)}\), uniformly over
the original flow and its restarts.
No uniform initial trace moment is required.
\end{lemma}

\begin{proof}
Apply Lemma~\ref{lem:weighted-energy} with \(w=Z\), \(\gamma=4\),
and
\[
 q=-c(H+\lambda_1)-\mathcal F
                  +C(1+|\partial g|+|\nabla\bar\nabla g|),
\]
as supplied by \eqref{eq:offcritical}. Discard the favorable
\(\mathcal F\) term to obtain \eqref{eq:largest-dissipation}.
The chain rule and the Cauchy--Schwarz inequality give
\[
\begin{aligned}
 \Gamma(Z^{p/2})
 &=\frac{p^2}{4}Z^p\Gamma(\log Z),\\
 \Gamma\!\left(e^{\psi_t^2/2}Z^{p/2}\right)
 &\le e^{\psi_t^2}Z^p
 \left[\frac{p^2}{2}\Gamma(\log Z)
       +2\psi_t^2\Gamma(\psi_t)\right].
\end{aligned}
\]
Since \(e^{\psi_t^2}\ge1\) and \(\psi_t\) is uniformly bounded by
Lemma~\ref{lem:flow-density}, these imply
\eqref{eq:eigenvalue-energy-control}.
\end{proof}

\begin{lemma}[Trace decay from a gradient modulus]\label{lem:trace-decay}
Let \(\psi\) be the original flow of Proposition~\ref{prop:common-flow}.
Suppose that, for some \(0\le\beta<1\) and \(K_\beta>0\),
\[
 \sup_{0\le t\le T_*}\|\psi(t)\|_{C^{1,\beta}(M)}\le K_\beta.
\]
Then, for each fixed \(\varepsilon>0\),
\begin{equation}\label{eq:trace-decay}
 \sup_M\tr A_\psi(t)\le C_{\beta,\varepsilon}t^{-\frac{1-\beta}{1+\beta}-\varepsilon},
 \qquad0<t\le T_*.
\end{equation}
The constant depends only on
\(M,\chi,\omega,n,k,\alpha,\|G\|_{C^\alpha(M)},\beta,\varepsilon,K_\beta\).
It is independent
of the individually finite initial moment.
\end{lemma}
\begin{proof}
For \(\psi_r=P_{r^2}\psi\), the assumed bound gives
\[
 \|\partial(\psi-\psi_r)\|_\infty\le C_\beta r^\beta,\qquad
 \|\Delta_\chi\psi_r\|_\infty\le C_\beta r^{\beta-1}.
\]
Fix the flow time and let \(0<r\le1\). Split
\(\tr A_\psi=\tr_\chi\omega+\Delta_\chi\psi_r+
\Delta_\chi(\psi-\psi_r)\) to obtain
\[
\begin{aligned}
 \int_M Z^p\tr A_\psi\,dV
 &=\int_M Z^p(\tr_\chi\omega+\Delta_\chi\psi_r)\,dV
   +\int_M Z^p\Delta_\chi(\psi-\psi_r)\,dV\\
 &\le Cr^{\beta-1} \int_M Z^p\,dV
   +\int_M Z^p\Delta_\chi(\psi-\psi_r)\,dV.
\end{aligned}
\]
Here the fixed background is bounded and \(r^{\beta-1}\ge1\).
By \eqref{eq:fixed-ibp} and the gradient approximation estimate,
\[
\begin{aligned}
 \int_M Z^p\Delta_\chi(\psi-\psi_r)\,dV
 &=-\int_M\Gamma_\chi(Z^p,\psi-\psi_r)\,dV\\
 &\le Cr^\beta\int_M|\partial Z^p|_\chi\,dV\\
 &=2Cr^\beta\int_M Z^{p/2}|\partial Z^{p/2}|_\chi\,dV.
\end{aligned}
\]
Lemma~\ref{lem:algebra} gives \(F^{i \bar i}\ge(\tr A_\psi)^{-1}\), hence
\[
 \frac{|\partial Z^{p/2}|_\chi^2}{\tr A_\psi}
 \le\Gamma(Z^{p/2}).
\]
Since \(d\nu=\Sk(A_\psi)\,dV\) and \(\Sk(A_\psi)\ge f_-\),
weighted Cauchy--Schwarz yields
\[
\begin{aligned}
 \int_M Z^{p/2}|\partial Z^{p/2}|_\chi\,dV
 &=\int_M(Z^p\tr A_\psi)^{1/2}
       \frac{|\partial Z^{p/2}|_\chi}{\sqrt{\tr A_\psi}}\,dV\\
 &\le (\int_M Z^p\tr A_\psi\,dV)^{1/2}
       \left(\int_M
       \frac{|\partial Z^{p/2}|_\chi^2}{\tr A_\psi}\,dV\right)^{1/2}\\
 &\le f_-^{-1/2}(\int_M Z^p\tr A_\psi\,dV)^{1/2}
       \left(\int_M\Gamma(Z^{p/2})\,d\nu\right)^{1/2}.
\end{aligned}
\]
Combining these estimates gives
\[
 \int_M Z^p\tr A_\psi\,dV\le Cr^{\beta-1}\int_M Z^p\,dV+
 Cr^\beta\left(\int\Gamma(Z^{p/2})\,d\nu\right)^{1/2}(\int_M Z^p\tr A_\psi\,dV)^{1/2}.
\]
The constants in this calculation depend only on
\(M,\chi,\omega,n,k,\beta,K_\beta,f_-\), independently of \(p,t,r\).
Let \(C_0\) be the constant in the preceding bound.
Since \(Z\asymp\tr A_\psi\), Jensen's inequality gives
\[
 \int_M Z^p\tr A_\psi\,dV \ge c\int_M Z^{p+1}\,dV
 \ge c\,\operatorname{Vol}(M,dV)^{-1/p}(\int_M Z^p\,dV)^{1+1/p}
 \ge c_p (\int_M Z^p\,dV)^{1+1/p}.
\]
If \( \int_M Z^p\tr A_\psi\,dV \ge2C_0 \int_M Z^p\,dV\), choose
\[
 r=\left(\frac{2C_0 \int_M Z^p\,dV}{\int_M Z^p\tr A_\psi\,dV}\right)^{1/(1-\beta)}\le1.
\]
Then \(C_0r^{\beta-1}\int_M Z^p\,dV=(\int_M Z^p\tr A_\psi\,dV)/2\). Moving this term to the left,
dividing by \((\int_M Z^p\tr A_\psi\,dV)^{1/2}\), and squaring gives
\[
\begin{aligned}
 \int_M\Gamma(Z^{p/2})\,d\nu
 &\ge\frac{\int_M Z^p\tr A_\psi\,dV}{4C_0^2r^{2\beta}}
 =c_\beta (\int_M Z^p\tr A_\psi\,dV)^{(1+\beta)/(1-\beta)}
       (\int_M Z^p\,dV)^{-2\beta/(1-\beta)}\\
 &\ge c_{p,\beta} (\int_M Z^p\,dV)^{1+(1+\beta)/(p(1-\beta))}.
\end{aligned}
\]
If \(\int_M Z^p\tr A_\psi\,dV<2C_0 \int_M Z^p\,dV\), the lower bound \(\int_M Z^p\tr A_\psi\,dV\ge c_p (\int_M Z^p\,dV)^{1+1/p}\)
instead gives \(\int_M Z^p\,dV \le(2C_0/c_p)^p\).
The nonnegativity of the energy therefore yields, in both cases,
\[
 \int_M\Gamma(Z^{p/2})\,d\nu
 \ge c_{p,\beta}(\int_M Z^p\,dV)^{1+(1+\beta)/(p(1-\beta))}-C_{p,\beta} \int_M Z^p\,dV.
\]

Fix \(p\ge p_0\) and write
\(\delta=(1+\beta)/(p(1-\beta))>0\).
The bounds for \(\psi_t\) and \(\Sk(A_\psi)\) give \(c \int_M Z^p\,dV \le Y_p\le C \int_M Z^p\,dV\).
Combining the last energy bound with
\eqref{eq:largest-dissipation} and
\eqref{eq:eigenvalue-energy-control}, and discarding the other
nonnegative dissipation terms, gives
\begin{equation}\label{eq:trace-moment-ode}
 Y_p'+c_{p,\beta}Y_p^{1+\delta}
 \le C_{p,\beta}(1+\rho(t))Y_p.
\end{equation}
Define the integrating factor by
\[
 \mathcal R(t)=C_{p,\beta}\int_0^t(1+\rho(s))\,ds.
\]
Since \(\rho(s)\le C(1+s^{\alpha/2-1})\) and \(\alpha>0\),
\(0\le\mathcal R(t)\le C_{p,\beta}(t+t^{\alpha/2})
\le C_{p,\beta}\) for \(0\le t\le T_*\le1\).
Multiplying \eqref{eq:trace-moment-ode} by \(e^{-\mathcal R}\)
and then differentiating the negative power yields
\[
\begin{aligned}
 \frac{d}{dt}(e^{-\mathcal R}Y_p)
 &\le-c_{p,\beta}e^{\delta\mathcal R}
             (e^{-\mathcal R}Y_p)^{1+\delta},\\
 \frac{d}{dt}(e^{-\mathcal R}Y_p)^{-\delta}
 &\ge c_{p,\beta}\delta e^{\delta\mathcal R}
 \ge c_{p,\beta}\delta.
\end{aligned}
\]
Integrating from \(0\) to \(t\) and using \(\mathcal R(0)=0\), we obtain
\begin{equation}\label{eq:trace-moment-decay}
\begin{aligned}
 Y_p(t)
 &\le e^{\mathcal R(t)}
       \left[Y_p(0)^{-\delta}+c_{p,\beta}\delta t\right]^{-1/\delta}\\
 &\le C_{p,\beta}t^{-1/\delta}
 =C_{p,\beta}t^{-p(1-\beta)/(1+\beta)}.
\end{aligned}
\end{equation}
Only finiteness of \(Y_p(0)\) is used: discarding
\(Y_p(0)^{-\delta}\ge0\) makes the final bound independent of
the initial moment. For \(\beta=0\), this gives \(Y_p(t)\le C_pt^{-p}\).

Combining (\ref{eq:largest-dissipation}) and (\ref{eq:eigenvalue-energy-control}), we can get:

\begin{equation*}
Y_p' + C \int_M [\Gamma(Z^{p/2})+ \Gamma(e^{\psi_t^2/2} Z^{p/2})]d\nu \le Cp \rho(t) Y_p
\end{equation*}

Using Proposition~\ref{prop:hessian-sobolev}, we can use the Sobolev inequality and apply the parabolic Moser Iteration to the above differential inequality to get that for any $1< \kappa < \frac{n+1}{n}$:
\[
 \|Z(t)\|_\infty\le
 C_p t^{-\kappa/((\kappa-1)p)}
 \left(\int_{t/2}^tY_p(s)\,ds\right)^{1/p}
 \le C_{p,\beta}t^{-\frac{1-\beta}{1+\beta}-1/((\kappa-1)p)}.
\]
In the second inequality above, we use (\ref{eq:trace-moment-decay}).
\end{proof}

\subsection{The gradient exponent bootstrap}

The following Lemma says that the blow up rate of $tr A_{\psi}(t)$ can imply the bound on the $C^{1,\beta'}$ norm of $\psi(t)$, which can be use the improve the blow up decay of $tr A_{\psi}(t)$.
\begin{lemma}[Trace improvement by dyadic restarts]\label{lem:dyadic-restart-improvement}
Let \(\psi\) and \(\varphi^\tau\) be the original flow and restarts
of Proposition~\ref{prop:common-flow}. Suppose that, for some
\(0<d<1+\alpha/2\) and \(C_{\mathrm{tr}}>0\),
\[
 \sup_M\tr A_\psi(t)\le C_{\mathrm{tr}}t^{-d},
 \qquad 0<t\le T_*.
\]
Then, for every \(0<\gamma<1\),
\begin{equation}\label{eq:restart-holder-from-trace}
 \sup_{0\le s\le T_*}\|\varphi^\tau(s)\|_{C^{1,\gamma}(M)}
 \le C_\gamma\tau^{-d},
 \qquad 0<\tau\le\tau_0.
\end{equation}
For every
\(0<\beta'<((1+\alpha/2)-d)/((1+\alpha/2)+d)\),
\begin{equation}\label{eq:dyadic-original-holder}
 \sup_{0\le s\le T_*}\|\psi(s)\|_{C^{1,\beta'}(M)}
 \le C_{\beta'}.
\end{equation}
Moreover, for every \(\varepsilon>0\),
\begin{equation}\label{eq:dyadic-trace-improvement}
 \sup_M\tr A_\psi(t)
 \le C_{d,\varepsilon}t^{-d/(1+\alpha/2)-\varepsilon},
 \qquad 0<t\le T_*.
\end{equation}
All constants depend only on
\(M,\chi,\omega,n,k,\alpha,\|G\|_{C^\alpha(M)},C_{\mathrm{tr}},d\)
and the indicated exponents \(\gamma,\beta',\varepsilon\).
\end{lemma}
\begin{proof}
By \eqref{eq:coarse-propagation} and
\(\varphi^\tau(0)=\psi(\tau)\),
\[
 \sup_{0\le s\le T_*}\sup_M\tr A_{\varphi^\tau}(s)
 \le C\sup_M\tr A_\psi(\tau)
 \le CC_{\mathrm{tr}}\tau^{-d}.
\]
The identity
\(\Delta_\chi\varphi^\tau=\tr A_{\varphi^\tau}-\tr_\chi\omega\),
the uniform potential bound, and fixed-metric Poisson estimates give
\[
 \|\varphi^\tau(s)\|_{C^{1,\gamma}}
 \le C_\gamma\bigl(\|\Delta_\chi\varphi^\tau(s)\|_\infty
                    +\|\varphi^\tau(s)\|_\infty\bigr)
 \le C_\gamma\tau^{-d}.
\]
This proves \eqref{eq:restart-holder-from-trace}.
The triangle inequality consequently gives, at the same elapsed time \(s\),
\[
\begin{aligned}
 \|\varphi^\tau(s)-\varphi^{\tau/2}(s)\|_{C^{1,\gamma}}
 &\le\|\varphi^\tau(s)\|_{C^{1,\gamma}}
       +\|\varphi^{\tau/2}(s)\|_{C^{1,\gamma}}\\
 &\le C_\gamma\bigl(\tau^{-d}+(\tau/2)^{-d}\bigr)
 \le C_\gamma\tau^{-d}.
\end{aligned}
\]
Similarly, \eqref{eq:restart-comparison} yields
\[
\begin{aligned}
 \|\varphi^\tau(s)-\varphi^{\tau/2}(s)\|_\infty
 &\le\|\varphi^\tau(s)-\psi(s)\|_\infty
       +\|\varphi^{\tau/2}(s)-\psi(s)\|_\infty\\
 &\le C\tau^{1+\alpha/2}+C(\tau/2)^{1+\alpha/2}
 \le C\tau^{1+\alpha/2}.
\end{aligned}
\]
Fix \(\beta'\) in the stated range and choose \(\beta'<\gamma<1\)
sufficiently close to \(1\) that
\((1+\alpha/2)(\gamma-\beta')-d(1+\beta')>0\).
Spatial interpolation gives
\[
\begin{aligned}
 \|\varphi^\tau(s)-\varphi^{\tau/2}(s)\|_{C^{1,\beta'}}
 &\le C_{\beta',\gamma}
   \|\varphi^\tau(s)-\varphi^{\tau/2}(s)\|_\infty^{(\gamma-\beta')/(1+\gamma)}\\
 &\qquad{}\times
   \|\varphi^\tau(s)-\varphi^{\tau/2}(s)\|_{C^{1,\gamma}}^{(1+\beta')/(1+\gamma)}\\
 &\le C_{\beta',\gamma}
 \tau^{[(1+\alpha/2)(\gamma-\beta')-d(1+\beta')]/(1+\gamma)}.
\end{aligned}
\]
For \(\tau_j=2^{-j}\tau_0\), the positive exponent makes these
consecutive differences summable in \(C^{1,\beta'}\), uniformly in \(s\).
The first restart is bounded by
\eqref{eq:restart-holder-from-trace}, and
\eqref{eq:restart-comparison} identifies the \(C^0\) limit as \(\psi(s)\).
This proves \eqref{eq:dyadic-original-holder}.

As \(\beta'\) approaches its upper endpoint,
\[
 \frac{1-\beta'}{1+\beta'}\longrightarrow
 \frac{d}{1+\alpha/2}.
\]
Given \(\varepsilon>0\), choose \(\beta'\) below that endpoint so that
\((1-\beta')/(1+\beta')<d/(1+\alpha/2)+\varepsilon/2\).
Applying Lemma~\ref{lem:trace-decay} with this \(\beta'\) and loss
\(\varepsilon/2\) proves \eqref{eq:dyadic-trace-improvement}.
\end{proof}

\begin{proposition}\label{lem:gradient-bootstrap}
For every fixed \(0<\beta<1\) and \(\varepsilon>0\), the original
flow on the interval of Proposition~\ref{prop:common-flow} satisfies
\begin{equation}\label{eq:all-beta}
 \sup_{0\le s\le T_*}\|\psi(s)\|_{C^{1,\beta}}\le C_\beta,
 \qquad \tr_\chi\omega_\psi(t)\le C_\varepsilon t^{-\varepsilon}.
\end{equation}
\end{proposition}
\begin{proof}
Proposition~\ref{prop:flow-bounds} bounds the gradient of the original
flow. Choose \(0<\varepsilon_0<\alpha/2\) and apply
Lemma~\ref{lem:trace-decay} with \(\beta=0\) to obtain
\[
 \sup_M\tr A_\psi(t)\le Ct^{-d_0},
 \qquad d_0=1+\varepsilon_0<1+\alpha/2.
\]
Fix \(q\) with \((1+\alpha/2)^{-1}<q<1\).

Apply Lemma~\ref{lem:dyadic-restart-improvement} by induction, we can prove the second inequality in (\ref{eq:all-beta}). Then, we apply the Lemma \ref{lem:dyadic-restart-improvement} to get the first inequality in (\ref{eq:all-beta}). 
\end{proof}

\section{A $W^{2,p}$ estimate for the initial potential}
\label{sec:finite}

\begin{lemma}\label{lem:coercivity}
For \(\mathsf H\in\Gamma_k\), \(k\ge2\), and every Hermitian \(B\),
\begin{equation}\label{eq:coercivity}
 -D^2\log\sigma_k(\mathsf H)[B,B]\ge |B|_{\mathrm F}^2/(\tr \mathsf H)^2.
\end{equation}
\end{lemma}
\begin{proof}
Hyperbolicity of \(\sigma_k\) with respect to \(\mathsf H\in\Gamma_k\) gives
real numbers \(\rho_1,\ldots,\rho_k\) such that
\[
 \sigma_k(\mathsf H+tB)=\sigma_k(\mathsf H)\prod_{\ell=1}^k(1+t\rho_\ell).
\]
Thus \(-D^2\log\sigma_k(\mathsf H)[B,B]=\sum\rho_\ell^2\).
If \(-D^2\log\sigma_k(\mathsf H)[B,B] r^2<1\), both \(\mathsf H+rB\) and \(\mathsf H-rB\) lie in the same
positive component and hence in \(\Gamma_2\). Consequently, we use Lemma \ref{lem:algebra} to get:
\[
 2r|B|_{\mathrm F}
 \le |\mathsf H+rB|_{\mathrm F}+|\mathsf H-rB|_{\mathrm F}< tr(\mathsf H+rB)+ tr(\mathsf H-rB) = 2\tr \mathsf H.
\]
Letting \(r\uparrow (-D^2\log\sigma_k(\mathsf H)[B,B])^{-1/2}\) proves \eqref{eq:coercivity}. If \(-D^2\log\sigma_k(\mathsf H)[B,B]=0\), arbitrary \(r\)
forces \(B=0\).
\end{proof}

\begin{lemma}[A superlevel estimate for the initial trace]\label{lem:initial-trace-superlevel}
Let \(2\le k\le n\), \(0<\alpha<1\), and let \(u\) be a smooth normalized
solution of \eqref{eq:main}. Let \(\psi\) be the original flow
\eqref{eq:original-flow} on the interval \([0,T_*]\) of
Proposition~\ref{prop:common-flow}. For every \(0<t\le T_*\),
\begin{equation}\label{eq:initial-trace-superlevel}
 \operatorname{Vol}_\chi\bigl\{x\in M:
       \tr A_u(x)>2\tr A_\psi(t,x)\bigr\}\le Ct^{\alpha/2}.
\end{equation}
Here \(\operatorname{Vol}_\chi(E)=\int_E\chi^n\), and \(C\) depends only on
\(M,\chi,\omega,n,k,\alpha,\|G\|_{C^\alpha(M)}\).
\end{lemma}
\begin{proof}
Fix \(0<t\le T_*\) and denote the set in \eqref{eq:initial-trace-superlevel}
by \(E_t\). The segment
\[
 H_\theta=(1-\theta)A_\psi(t)+\theta A_u,\qquad 0\le\theta\le1,
\]
lies in the G{\aa}rding cone by its convexity
\cite[Theorem~2.2]{Spruck2005}. Since \(F=\log\Sk\) differs from
\(\log\sigma_k\) by a constant, Lemma~\ref{lem:coercivity} gives
\[
 -\frac{d^2}{d\theta^2}F(H_\theta)
 \ge \frac{|A_u-A_\psi(t)|_{\mathrm F}^2}{(\tr H_\theta)^2}
 \ge 0.
\]
On \(E_t\), we have \(\tr A_\psi(t)<\tfrac12\tr A_u\) and
\(\tr A_u>0\) by admissibility. Linearity of the trace and
\(0\le\theta\le1\) give
\[
 \tr H_\theta
 =(1-\theta)\tr A_\psi(t)+\theta\tr A_u
 \le \left(\frac{1-\theta}{2}+\theta\right)\tr A_u
 =\frac{1+\theta}{2}\tr A_u\le\tr A_u.
\]
Also, Cauchy--Schwarz in a \(\chi\)-unitary frame gives
\[
 |A_u-A_\psi(t)|_{\mathrm F}
 \ge \frac{\tr A_u-\tr A_\psi(t)}{\sqrt n}
 \ge \frac{\tr A_u}{2\sqrt n}\qquad\text{on }E_t.
\]
Thus \(-\frac{d^2}{d\theta^2}F(H_\theta)\ge\mathbf{1}_{E_t}/(4n)\)
on \(M\).

Let \(L_t\) be the linearization at \(A_\psi(t)\). Since
\(A_u-A_\psi(t)=\chi^{-1}\ddc(u-\psi(t))\), the derivative of
\(F(H_\theta)\) at \(\theta=0\) is \(L_t(u-\psi(t))\).
Taylor's formula and \(F(A_u)=G\) now give
\[
 F(A_\psi(t))+L_t(u-\psi(t))-G
 =-\int_0^1(1-\theta)\frac{d^2}{d\theta^2}F(H_\theta)\,d\theta
 \ge \frac{\mathbf{1}_{E_t}}{8n}.
\]
Integrate against \(d\nu_t=\Sk(A_\psi(t))\chi^n\). Since
\(\omega_\psi(t)\) is closed, Lemma~\ref{lem:newton-symmetry} gives
\(\int_M L_t(u-\psi(t))\,d\nu_t=0\). By \eqref{eq:original-flow},
\[
 \frac{\nu_t(E_t)}{8n}
 \le \int_M\bigl(F(A_\psi(t))-G\bigr)\,d\nu_t
 =\int_M\bigl(\psi_t(t)+P_tG-G\bigr)\,d\nu_t.
\]
Proposition~\ref{prop:flow-bounds} and Lemma~\ref{lem:heat-regularization}
give \(\|\psi_t(t)\|_\infty+\|P_tG-G\|_\infty\le Ct^{\alpha/2}\).
Using \(f_-\chi^n\le d\nu_t\le f_+\chi^n\) from
\eqref{eq:common-flow-uniform-bounds} proves \eqref{eq:initial-trace-superlevel}.
\end{proof}

\begin{proposition}[Finite initial Sobolev bounds]\label{prop:finite-moments}
Let \(0<\alpha<1\). Every smooth normalized solution \(u\) of \eqref{eq:main} satisfies
\begin{equation}\label{eq:finite-moments}
 \|u\|_{W^{2,p}}\le C_p,\qquad1\le p<\infty.
\end{equation}
Here \(C_p\) depends only on
\(M,\chi,\omega,n,k,\alpha,p,\|G\|_{C^\alpha(M)}\).
\end{proposition}
\begin{proof}
For \(k=1\), \eqref{eq:main} gives
\(\Delta_\chi u=ne^G-\tr_\chi\omega\), and the Schauder estimate for
the normalized Poisson equation gives \eqref{eq:finite-moments}.
Suppose \(2\le k\le n\). Fix \(2\le p<\infty\) and choose
\(\varepsilon=\alpha/(4p)\). Proposition~\ref{lem:gradient-bootstrap}
gives \(\sup_M\tr A_\psi(t)\le C_\varepsilon t^{-\varepsilon}\)
for \(0<t\le T_*\). For every
\(\lambda\ge\lambda_0:=2C_\varepsilon T_*^{-\varepsilon}\), take
\[
 t_\lambda=\left(\frac{2C_\varepsilon}{\lambda}\right)^{1/\varepsilon}.
\]
Then \(0<t_\lambda\le T_*\) and \(2\tr A_\psi(t_\lambda)\le\lambda\)
everywhere on \(M\). Thus Lemma~\ref{lem:initial-trace-superlevel} gives
\[
\begin{aligned}
 \operatorname{Vol}_\chi\{\tr A_u>\lambda\}
 &\le\operatorname{Vol}_\chi\{\tr A_u>2\tr A_\psi(t_\lambda)\}\\
 &\le C t_\lambda^{\alpha/2}\le C_p\lambda^{-2p},
\end{aligned}
\]
where the last inequality uses \(\alpha/(2\varepsilon)=2p\).
Since \(\tr A_u>0\), the layer-cake formula yields
\[
\begin{aligned}
 \int_M(\tr A_u)^p\,\chi^n
 &=p\int_0^\infty \lambda^{p-1}
        \operatorname{Vol}_\chi\{\tr A_u>\lambda\}\,d\lambda\\
 &\le \operatorname{Vol}_\chi(M)\lambda_0^p
     +pC_p\int_{\lambda_0}^\infty\lambda^{-p-1}\,d\lambda
 \le C_p.
\end{aligned}
\]
Finally, \(\Delta_\chi u=\tr A_u-\tr_\chi\omega\), and
Proposition~\ref{prop:flow-bounds} at \(t=0\) bounds \(\|u\|_\infty\).
The elliptic \(L^p\) estimate for \(\Delta_\chi\) gives
\[
 \|u\|_{W^{2,p}(M,\chi)}
 \le C_{M,\chi,p}\bigl(\|\Delta_\chi u\|_{L^p(\chi^n)}
                       +\|u\|_{L^p(\chi^n)}\bigr)\le C_p.
\]
The range \(1\le p<2\) follows from the case \(p=2\) by H\"older's
inequality on \(M\).
\end{proof}

\section{$W^{2,p}$ estimate and trace propagation along flow}\label{sec:additive}
We return to the original flow \(\psi\) and its restarts
\(\varphi^\tau\), with endomorphisms \(A_\psi(t)\) and
\(A_{\varphi^\tau}(s)\) from Definition~\ref{def:flows}.
Their measures are
\[
 d\nu_t=\Sk(A_\psi(t))\chi^n,\qquad
 d\nu_s^\tau=\Sk(A_{\varphi^\tau}(s))\chi^n.
\]
Proofs are written for \(\psi\); on a restart, replace
\(\psi(t),A_\psi(t),\psi_t(t),\nu_t\) by
\(\varphi^\tau(s),A_{\varphi^\tau}(s),
\partial_s\varphi^\tau(s),\nu_s^\tau\).
The operators \(a,L,\Gamma,D\) are evaluated on the same path.
The estimates of this section use the initial Sobolev bounds only
after Proposition~\ref{prop:finite-moments} has established them.

\begin{lemma}[Uniform trace moments]\label{lem:restart-moments}
The original flow and all restarts
satisfy
\begin{equation}\label{eq:restart-moments}
\begin{aligned}
 \sup_{0\le t\le T_*}\|\tr A_\psi(t)\|_{L^p(\nu_t)}
 &\le C_p,\\
 \sup_{0<\tau\le\tau_0}\sup_{0\le s\le T_*}
 \|\tr A_{\varphi^\tau}(s)\|_{L^p(\nu_s^\tau)}
 &\le C_p,\qquad1\le p<\infty.
\end{aligned}
\end{equation}
\end{lemma}
\begin{proof}
Only now do we use the initial moments
\eqref{eq:finite-moments}. Dropping favorable terms in
\eqref{eq:largest-dissipation} gives
\(Y_p'\le Cp\rho Y_p\), with integrable \(\rho\).
Since \(Z\asymp\tr A_\psi\) and all density and velocity weights are
bounded, Gronwall first yields
\(\sup_{t\le T_*}\|\operatorname{tr}_{\chi}\omega_{\psi}(t)\|_{L^p}\le C_p\).
Every restart begins at an original time
\(\tau\le\tau_0\le T_*/2\). Its only different initial entropy
weight is the uniformly bounded, generally nonzero velocity.
The same Gronwall argument in the restart's own clock gives
\eqref{eq:restart-moments}.

Lower exponents follow by H\"older. No limiting statement as
\(p\to\infty\) is inferred.
\end{proof}

\begin{lemma}[Trace subsolution]\label{lem:trace-subsolution}
On the original flow,
\begin{equation}\label{eq:trace-subsolution}
 D(\tr A_\psi)\le C(\mathcal F)\tr A_\psi+|\Delta g|
       \le C(\tr A_\psi)^k+|\Delta g|.
\end{equation}
The same inequalities hold on a restart with \(A_\psi\) replaced by
\(A_{\varphi^\tau}\), \(g=P_sG\), and the restart operators.
The sources \((\tr A_\psi)^k\) and
\((\tr A_{\varphi^\tau})^k\) have every fixed finite moment uniformly.
\end{lemma}
\begin{proof}
Since \(\Omega=\omega+\ddc\psi\) is closed, K\"ahler commutation gives
\(\Omega_{i\bar i;p\bar p}-\Omega_{p\bar p;i\bar i}=O(|A_\psi|)\).
Differentiating \(\psi_t=F(A_\psi)-g\) and tracing yields
\[
 D(\tr A_\psi)=\sum_pF''[\nabla_p A_\psi,\nabla_{\bar p} A_\psi]
              +O((\mathcal F)|A_\psi|)-\Delta g.
\]
Concavity discards the first term. Lemma~\ref{lem:algebra} gives
\(|A_\psi|\le\tr A_\psi\), and \eqref{eq:trace-growth} gives
\(\mathcal F\le C(\tr A_\psi)^{k-1}\). This proves
\eqref{eq:trace-subsolution}; Lemma~\ref{lem:restart-moments}
supplies the source moments.
\end{proof}

\begin{proposition}[Forward heat smoothing]\label{prop:forward-smoothing}
Along the original flow \(\psi_t=F(A_\psi)-P_tG\) of
\eqref{eq:original-flow}, take
\[
 L_t h=F^{i\bar j}(A_\psi(t))\nabla_i\nabla_{\bar j}h,\qquad
 d\nu_t=\Sk(A_\psi(t))\chi^n .
\]
On an unshifted restart, replace \(A_\psi(t)\) and \(\nu_t\) throughout
by \(A_{\varphi^\tau}(t)\) and \(\nu_t^\tau\), using the restart's
own time variable. For \(h\in C(M)\), let \(w(t,\cdot)=T(t,s)h\) solve
\[
 \begin{cases}
  (\partial_t-L_t)w=0,&s<t\le T_*,\\
  w(s,\cdot)=h.
 \end{cases}
\]
In particular \(T(s,s)h=h\).

Let \(B\) be a common upper bound for \(|\psi_t|\) and
\(|\partial_t\varphi^\tau|\) on \(M\times[0,T_*]\), as supplied by
Lemma~\ref{lem:flow-density}. For \(1<\sigma<n/(n-1)\) and
\(p_0\le p<\infty\), the operator extends to \(L^p(\nu_s)\) and satisfies
\begin{equation}\label{eq:forward-smoothing}
 \|T(t,s)h\|_\infty\le C_p(t-s)^{-\sigma/((\sigma-1)p)}
                  \|h\|_{L^p(\nu_s)},\qquad 0\le s<t\le T_* .
\end{equation}
One may take \(p_0=4B^2+2B+2\). By \eqref{eq:flow-velocity-bounds},
\(B\) can be chosen in terms of \(M,\chi,n,\alpha,\|G\|_{C^\alpha(M)}\).
The constant \(C_p\) depends only on
\(M,\chi,\omega,n,k,\alpha,\|G\|_{C^\alpha(M)},\sigma,p\), uniformly
over these paths. For every \(0\le s\le t\le T_*\), \(T(t,s)\) is positive and
\begin{equation}\label{eq:markov}
 T(t,s)1=1,\qquad \|T(t,s)h\|_\infty\le\|h\|_\infty
 \quad\text{for every }h\in C(M).
\end{equation}
\end{proposition}

\begin{proof}
On each smooth admissible path the coefficients of \(L_t\) are smooth
and uniformly elliptic on \(M\times[0,T_*]\), with ellipticity
constants allowed to depend on the path. The scalar maximum principle
gives positivity and \eqref{eq:markov}, since \(L_t1=0\).

Consider the original flow and take positive smooth initial data \(h\).
Then \(w=T(t,s)h>0\) satisfies \(D\log w=\Gamma(\log w)\).
Lemma~\ref{lem:weighted-energy} applies with \(q=0\) and \(\gamma=1\),
which gives the stated choice of \(p_0\). The chain rule gives
\[
 \Gamma(e^{\psi_t^2/2}w^{p/2})
 \le e^{\psi_t^2}w^p
  \left(\tfrac12 p^2\Gamma(\log w)+2B^2\Gamma(\psi_t)\right).
\]
Hence, for \(Y_p(t)=\int_M e^{\psi_t^2}w^p\,d\nu_t\),
\begin{equation}\label{eq:forward-energy}
 Y_p'+c\int_M\Gamma(e^{\psi_t^2/2}w^{p/2})\,d\nu_t
 \le C(1+\|\partial_tP_tG\|_\infty)Y_p,\qquad p\ge p_0 .
\end{equation}
Here \(c,C\) depend only on \(B\). By Lemma~\ref{lem:heat-regularization},
\(\int_0^{T_*}(1+\|\partial_rP_rG\|_\infty)\,dr\) is uniformly bounded,
so an integrating factor absorbs the right-hand side. The spatial
Sobolev inequality of Proposition~\ref{prop:hessian-sobolev}, with
exponent \(2\sigma\), and \(1\le e^{\psi_t^2}\le e^{B^2}\) give the
time-cutoff Moser estimate as in \cite[Chapter~19]{PeterLi2012},
\[
 \|w(t)\|_\infty
 \le C_p(t-s)^{-(2\sigma-1)/((\sigma-1)p)}
       \left(\int_s^t\int_M w(r)^p\,d\nu_r\,dr\right)^{1/p}.
\]
Gronwall's inequality in \eqref{eq:forward-energy} bounds the spacetime
integral by \(C(t-s)\|h\|_{L^p(\nu_s)}^p\), proving
\eqref{eq:forward-smoothing}. Positivity gives
\(|T(t,s)h|\le T(t,s)|h|\), and approximation extends the estimate to
all \(h\in L^p(\nu_s)\). The same proof uses
\(\partial_t\varphi^\tau\) and \(d\nu_t^\tau\) on each restart.
\end{proof}

\begin{proposition}[Additive trace propagation]\label{prop:additive}
Uniformly on the original flow and every unshifted restart,
\begin{equation}\label{eq:additive}
\begin{aligned}
 \sup_M\tr A_\psi(t)&\le\sup_M\tr A_\psi(s)+C(t-s)^{\alpha/2},\\
 \sup_M\tr A_{\varphi^\tau}(t)&\le
 \sup_M\tr A_{\varphi^\tau}(s)+C(t-s)^{\alpha/2},
 \qquad0\le s<t\le T_*.
\end{aligned}
\end{equation}
The second line uses the restart's own time.
\end{proposition}
\begin{proof}
We prove the estimate for \(\psi\); the same argument on
\(\varphi^\tau\) gives the second line.
For \(0<s<t\le T_*\), Duhamel's formula and comparison applied
to \eqref{eq:trace-subsolution} give the pointwise inequality
\[
\begin{aligned}
 \tr A_\psi(t)
 &\le T(t,s)[\tr A_\psi(s)]
       +C\int_s^t T(t,r)[(\tr A_\psi(r))^k]\,dr\\
 &\quad+\int_s^t T(t,r)[|\Delta g(r)|]\,dr.
\end{aligned}
\]
Indeed, the right-hand side solves the linear equation with source
\(C(\tr A_\psi)^k+|\Delta g|\) and initial value \(\tr A_\psi(s)\),
whereas \(\tr A_\psi\) is a subsolution. Positivity of the evolution
operator follows from the maximum principle.

Fix \(2<2\sigma<2n/(n-1)\) and choose one finite exponent
\[
 p>\max\left\{p_0,
       \frac{\sigma}{(\sigma-1)(1-\alpha/2)}\right\}.
\]
The initial term is bounded by \eqref{eq:markov}:
\[
 \sup_M T(t,s)[\tr A_\psi(s)]\le\sup_M\tr A_\psi(s).
\]
For the nonlinear source, \eqref{eq:forward-smoothing} and
Lemma~\ref{lem:restart-moments} at exponent \(kp\) give
\[
\begin{aligned}
 \|T(t,r)[(\tr A_\psi(r))^k]\|_\infty
 &\le C_p(t-r)^{-\sigma/((\sigma-1)p)}
            \|(\tr A_\psi(r))^k\|_{L^p(\nu_r)}\\
 &=C_p(t-r)^{-\sigma/((\sigma-1)p)}
            \|\tr A_\psi(r)\|_{L^{kp}(\nu_r)}^k\\
 &\le C_p(t-r)^{-\sigma/((\sigma-1)p)}.
\end{aligned}
\]
Since \(\sigma/((\sigma-1)p)<1-\alpha/2<1\), this is integrable,
and hence
\[
\begin{aligned}
 C\int_s^t\|T(t,r)[(\tr A_\psi(r))^k]\|_\infty\,dr
 &\le C_p\int_s^t(t-r)^{-\sigma/((\sigma-1)p)}\,dr\\
 &\le C_p(t-s)^{1-\sigma/((\sigma-1)p)}.
\end{aligned}
\]
For the remaining source, \(g(r)=P_rG\) in the path's own
time variable. Thus \eqref{eq:markov} and \eqref{eq:heat-envelope}
give
\[
 \int_s^t\|T(t,r)[|\Delta g(r)|]\|_\infty\,dr
 \le\int_s^t\|\Delta g(r)\|_\infty\,dr
 \le C\int_s^t r^{\alpha/2-1}\,dr.
\]
Combining these three estimates yields
\begin{align*}
 \sup_M\tr A_\psi(t)
 \le\sup_M\tr A_\psi(s)+C_p(t-s)^{1-\sigma/((\sigma-1)p)}
                   +C\int_s^t r^{\alpha/2-1}\,dr.
\end{align*}
Our choice of \(p\) gives
\(1-\sigma/((\sigma-1)p)>\alpha/2\), so the middle term is
at most \(C(t-s)^{\alpha/2}\) on the fixed time interval.
Moreover, subadditivity of \(r\mapsto r^{\alpha/2}\) gives
\[
 \int_s^t r^{\alpha/2-1}\,dr
 =\frac2\alpha(t^{\alpha/2}-s^{\alpha/2})
 \le\frac2\alpha(t-s)^{\alpha/2}.
\]
This proves \eqref{eq:additive}. The exponent \(p\) is fixed
once and for all, and its dependence is absorbed into \(C\);
no limit \(p\to\infty\) is used. The coefficient of the initial
supremum is exactly one because of \eqref{eq:markov}.
All formulas retain the actual possibly nonzero initial
restart velocity. Closed endpoints follow by continuity
and the integrability of \eqref{eq:heat-envelope}.
\end{proof}

\section{A uniform bound for \texorpdfstring{$\tr_\chi\omega_u$}{the initial trace}}\label{sec:initial-trace}

Let \(\psi\) be the original flow \eqref{eq:original-flow} with
\(\psi(0)=u\), and let \(m(s)=\max_M\tr A_{\psi(s)}\) for \(0\le s\le T_*\).
Write \(D_s=\partial_s-L_s\), where \(L_s\) is the linearization of
\(F=\log\Sk\) at \(A_{\psi(s)}\), and write \(r_+=\max\{r,0\}\).
We compare the flow on \([t/2,t]\) with the fixed potential \(\psi(2t)\).

\begin{lemma}[Comparison with the potential at time \(2t\)]\label{lem:trace-time-comparison}
Let \(2\le k\le n\) and \(0<t\le T_*/2\). For \(t/2\le s\le t\),
\begin{align}
 \|\psi(s)-\psi(2t)\|_\infty
 &\le Ct^{1+\alpha/2},\label{eq:trace-frozen-potential}\\
 D_s\bigl(\psi(s)-\psi(2t)\bigr)
 &\ge \frac{(\tr A_{\psi(s)}-m(2t))_+^2}
              {2n(\tr A_{\psi(s)})^2}-Ct^{\alpha/2}.
 \label{eq:trace-frozen-concavity}
\end{align}
The constant depends only on \(M,\chi,\omega,n,k,\alpha,\|G\|_{C^\alpha(M)}\).
\end{lemma}

\begin{proof}
Integrating the velocity bound \(\|\psi_t(r)\|_\infty\le Cr^{\alpha/2}\)
from Lemma~\ref{lem:flow-density} over \([s,2t]\) gives
\eqref{eq:trace-frozen-potential}. 

At each \(x\in M\), the segment
\(H_\theta=(1-\theta)A_{\psi(s)}+\theta A_{\psi(2t)}\),
\(0\le\theta\le1\), lies in the convex cone \(\Gamma_k\). The derivative
of \(F(H_\theta)\) at \(\theta=0\) equals \(L_s(\psi(2t)-\psi(s))\).
Holding \(t\) fixed, we have that $\partial_s \psi(2t)=0$. Taylor's formula and the flow equation at \(s\) give
\begin{equation}\label{e Dspsi}
 D_s(\psi(s)-\psi(2t))
 =F(A_{\psi(2t)})-P_sG
  -\int_0^1(1-\theta)\frac{d^2}{d\theta^2}F(H_\theta)\,d\theta.
\end{equation}
At time \(2t\) the flow equation gives
\(F(A_{\psi(2t)})=\psi_t(2t)+P_{2t}G\). Hence \eqref{eq:heat-envelope} yields
\[
 \|F(A_{\psi(2t)})-P_sG\|_\infty
 \le \|\psi_t(2t)\|_\infty+\|P_{2t}G-G\|_\infty
                         +\|P_sG-G\|_\infty\le Ct^{\alpha/2}.
\]
This estimates the first and the second term on the right-hand side of (\ref{e Dspsi}).

\textbf{Case 1}:  \(\tr A_{\psi(s)}>m(2t)\), we have
\(\tr A_{\psi(2t)}\le m(2t)<\tr A_{\psi(s)}\). Thus
\[
 \tr H_\theta\le\tr A_{\psi(s)},\qquad
 |A_{\psi(s)}-A_{\psi(2t)}|_{\mathrm F}
 \ge \frac{\tr A_{\psi(s)}-m(2t)}{\sqrt n}.
\]
Since \(F=\log\sigma_k-\log\binom nk\), Lemma~\ref{lem:coercivity} gives
\[
 -\frac{d^2}{d\theta^2}F(H_\theta)
 \ge\frac{|A_{\psi(s)}-A_{\psi(2t)}|_{\mathrm F}^2}{(\tr H_\theta)^2}.
\]
Integration against \(1-\theta\) gives the factor \(1/2\) in
\eqref{eq:trace-frozen-concavity}.

\textbf{Case 2}: \(\tr A_{\psi(s)}\le m(2t)\),
the positive part vanishes and concavity gives the same bound.
\end{proof}

\begin{proposition}[Initial trace bound]\label{prop:trace-cap}
Every smooth normalized admissible solution \(u\) of \eqref{eq:main} satisfies
\[
 \max_M\tr_\chi\omega_u\le C(M,\chi,\omega,n,k,\alpha,\|G\|_{C^\alpha(M)}).
\]
\end{proposition}

\begin{proof}
For \(k=1\), \eqref{eq:main} gives \(\tr_\chi\omega_u=ne^G\).
Suppose \(k\ge2\) and fix
\(\varepsilon=\min\{\alpha/12,1/[2(k-1)]\}\).
By Proposition \ref{lem:gradient-bootstrap}, we have
\[
 m(s)\le  C_\varepsilon s^{-\varepsilon},\qquad 0<s\le T_* .
\]

Fix \(0<t\le T_*/2\). For \(s\in [t/2,t]\), choose
\(K=C_\varepsilon(t/2)^{-\varepsilon}\ge1\) so that \(m(s)\le K\).
Then \(tK^{k-1}\le C\) by the choice of \(\varepsilon\). Put
\(E(s,x)=(\tr A_{\psi(s)}(x)-m(2t))_+\). We apply the maximum principle to
\[
 Q(s,x)=(s-t/2)E(s,x)^3-C_0K^3\bigl(\psi(s,x)-\psi(2t,x)\bigr).
\]
The function \(r\mapsto(r_+)^3\) is \(C^2\) and convex. Its derivative
gives the factor \(3E^2\), so Lemma~\ref{lem:trace-subsolution} implies
\[
\begin{aligned}
 D_s\bigl((s-t/2)E^3\bigr)
 &\le E^3+3(s-t/2)E^2D_s(\tr A_{\psi(s)})\\
 &\le E^2\bigl(K+CtK^k+Ct^{\alpha/2}\bigr)\le CKE^2 .
\end{aligned}
\]
By \eqref{eq:trace-frozen-concavity},
\(D_s(\psi(s)-\psi(2t))\ge E^2/(2nK^2)-Ct^{\alpha/2}\).
Take \(C_0>2nC\), with \(C\) from the preceding display. We can use the above estimate to get:
\begin{equation*}
D_s Q \le CK E^2 -C_0 K^3(\frac{E^2}{2nK^2}- C t^{\frac{\alpha}{2}})\le C C_0 K^3 t^{\alpha/2}.
\end{equation*} 
The maximum principle yields
\[
 \max_M Q(t)\le \max_M Q(t/2)+CK^3t^{1+\alpha/2}.
\]
The cubic term $(s-t/2)E(s,x)^3$ vanishes at \(t/2\), and \eqref{eq:trace-frozen-potential}
bounds the potential difference at both endpoints. Consequently
\[
 \frac t2\bigl(m(t)-m(2t)\bigr)_+^3\le CK^3t^{1+\alpha/2}.
\]
Taking cube roots gives \(m(t)\le m(2t)+CKt^{\alpha/6}\).
Since \(K=C_\varepsilon(t/2)^{-\varepsilon}\) and \(\varepsilon\le\alpha/12\),
\begin{equation}\label{eq:trace-two-times}
 m(t)\le m(2t)+Ct^{\alpha/12},\qquad 0<t\le T_*/2 .
\end{equation}
For \(t_j=2^{-j}T_*/2\), summing \eqref{eq:trace-two-times} gives
\(m(t_j)\le m(T_*/2)+C\sum_{\ell=1}^j t_\ell^{\alpha/12}\le C\).
The positive-time estimate above bounds \(m(T_*/2)\), since \(T_*\)
depends only on the quantities listed in the statement.
Smoothness at \(t=0\) and \(\psi(0)=u\) give
\(m(t_j)\to\max_M\tr_\chi\omega_u\).
\end{proof}

\section{Proof of Theorem \ref{thm:main}}\label{sec:sharp}

\begin{proof}[Proof of Theorem~\ref{thm:main}]
 Let $G_j$ be a sequence of smooth functions converging uniformly to $G$ with $||G_j||_{C^{\alpha}}$ uniformly bounded by a constant depending on $||G||_{C^{\alpha}}$ and $\int_M e^{G_j} \chi^n = \int_M \omega^k \wedge \chi^{n-k}$. Then, we can use the classical results by Dinew-Kolodziej \cite{DK2014} to obtain smooth solution $u_j$ solving 
\begin{equation*}
\omega_{u_j}^k \wedge \chi^{n-k}= e^{G_j}\chi^n, \lambda_{\chi} (\omega_{u_j})\in \Gamma_k, \sup_M u_j =0.
\end{equation*}
Using Proposition \ref{prop:trace-cap}, there exists a constant $C$ such that:
\begin{equation*}
tr_{\chi} \omega_{u_j}\le C
\end{equation*}
for any $j$. Using the Theorem 1.1 in \cite{Chu}, there exists a constant $C_1$ such that:
\begin{equation*}
||u_j||_{C^{2,\alpha}}\le C_1
\end{equation*}
for any $j$. Then, we can take a subsequence of $u_j$, converging to a function $u\in C^{2,\alpha}$ in $C^{2,\alpha'}$ for some $\alpha' < \alpha$.  $u$ solves (\ref{eq:main}). The uniqueness of the solution follows by maximum principle.
\end{proof}

\bibliographystyle{amsplain}
\bibliography{reference}

\begin{center}
\begin{minipage}{0.88\textwidth}
\small
\textbf{Yulun Xu}\\[0.35em]
University of Toronto\\
 40 St. George Street, Toronto, ON, Canada\\[0.35em]
\href{mailto:yulun.1.xu@gmail.com}{\texttt{yulun.1.xu@gmail.com}}
\end{minipage}
\end{center}
\end{document}